\documentclass[10pt,twoside]{article}
\usepackage{amssymb,amsmath}
\usepackage{amsfonts}

\usepackage[utf8]{inputenc}

\usepackage{makeidx}

\numberwithin{equation}{section}

\begin{document}

\title{Non-central limit theorem for non-linear functionals
of vector valued Gaussian stationary random fields}

\author{P\'eter Major \\
Alfr\'ed R\'enyi Institute of Mathematics \\
Budapest, P.O.B. 127 H--1364, Hungary, e-mail:major@renyi.hu}

\maketitle 

\medskip\noindent
{\it Abstract.}\/ Here I prove non-central limit theorems for
non-linear functionals of vector valued stationary random fields
under appropriate conditions. They are the multivariate versions 
of the results in paper~\cite{6}. Previously  A. M. Arcones 
formulated such a result in Theorem~6 of his paper~\cite{1}. 
But there are serious problems with his result. Even its
formulation must be corrected. I explain the problems related
to Arcones' paper in the main text. In this paper I present
the right formulation of the multivariate version of the
non-central limit theorem in paper~\cite{6} together with
its correct proof. To do this first the theory of the Gaussian
stationary random fields described in the work~\cite{9} had to
be generalized to the case of vector valued random fields. This
was done in my work published in two subsequent papers~\cite{10}
and~\cite{11}. Here I prove the multivariate version of the
result about non-central limit theorems in paper~\cite{6} with
their help.

\section{On the motivation for this research.}

In this paper the following problem is considered.

\medskip
Let us have a $d$-dimensional vector valued Gaussian stationary
random field $X(p)=(X_1(p),\dots,X_d(p))$, $p\in{\mathbb Z}^\nu$,
where ${\mathbb Z}^\nu$ denotes the  lattice points with integer
coordinates in the $\nu$-dimensional Euclidean space
${\mathbb R}^\nu$ and a function $H(x_1,\dots,x_d)$
of $d$ variables with arguments $x_s\in{\mathbb R}^\nu$,
$1\le s\le d$. We define with their help the random variables 
$Y(p)=H(X_1(p),\dots,X_d(p))$ for all $p\in{\mathbb Z}^\nu$.
Let us introduce for all $N=1,2,\dots$ the normalized sum
$$
S_N=A_N^{-1}\sum_{p\in B_N} Y(p) 
$$
with an appropriate norming constant $A_N>0$, where
\begin{equation}
 B_N=\{p=(p_1,\dots,p_\nu)\in{\mathbb Z}^\nu\colon\; 0< p_k\le N
\textrm { for all } 1\le k\le\nu\}. \label{1.1}
\end{equation}

In this paper a non-Gaussian limit theorem is proved for these 
normalized sums $S_N$ with an appropriate norming constant 
$A_N$ if this vector valued Gaussian stationary random field 
$X(p)$, $p\in{\mathbb Z}^\nu$, and the function $H(x_1,\dots,x_d)$  
satisfy certain conditions. Paper~\cite{6} contains such limit 
theorems for non-linear functionals of scalar valued stationary 
Gaussian random fields, and here their natural multivariate
generalizations are presented.

A. M. Arcones formulated such a result in Theorem~6 of
paper~\cite{1}, but I found his discussion unsatisfactory.
Here I explain the main problems related to his proof.

In the proof of the limit theorem in Theorem~6 of~\cite{1}
the spectral representation of the covariance function of a
vector valued stationary random process is needed. This
representation is presented in formula~(3.2) of~\cite{1}.
But the properties of the measures (actually complex
measures) of $G^{(p,q)}$ are not discussed. The same can be
told about the random spectral measures $Z_{G^{(p,p)}}$ in the
next formula~(3.3). These objects were defined in the
scalar valued case, and their basic properties were also
proved. But the generalization of these definitions to
the vector valued case and the proof of their properties
are far from trivial.

The same can be said about the statements of paper~\cite{1}
in formulas (3.5), (3.6) and (3.7). Here first the statements
must be corrected. In formula~(3.5) the arguments
$A_i$, $1\le i\le d$,  must be replaced by $\frac{A_i}n$.
Then it must be explained what kind of limit is taken
in this formula. Finally, it must also be explained what
kind of limit $(Z_{G^{(1,1)}_0},\dots,Z_{G^{(d,d)}_0})$ random
spectral measures appear here as the limit. It is not a
random spectral measure in the classical sense, it can
be interpreted only as the random spectral measure of a
generalized random field. The necessary definitions and
proofs are missing again.

Formulas~(3.6) and~(3.7) in~\cite{1} contain a limit theorem
which is actually a special case of Theorem~6. The main step
in Arcones' proof consists in the reduction of the result
in Theorem~6 to this special case. But the proof of the
result formulated in~(3.6) and~(3.7) is missing. This is
a limit theorem for a sequence of random vectors. The
convergence of the single coordinates of these vectors
follows from the already proved result in~\cite{6} which
deals with the one-dimensional version of this problem.
(More precisely, this one-dimensional convergence would
follow from this already proved result if formula~(3.7)
were written in the correct form. The random integrals
defined in it should be taken on $\mathbb R^\tau$ instead
of $[-\pi,\pi]^\tau$.) But I do not see how the results
proven in the one-dimensional case could help in the proof
the convergence of the random vectors, i.e. how the
result formulated in~(3.6) and~(3.7) could be proved by the
methods of~\cite{1}. I have the impression that the
proof of these formulas is not simpler than a direct
proof of Theorem~6 in~\cite{1}.

An appropriate proof of the non-central limit theorem
should start with a good and complete description of the
spectral representation of the covariance function of vector
valued stationary processes. This is done e.g. in
paper~\cite{5} of Cramer or in paper~\cite{13} of Rozanov.
This result is missing from Arcones' paper.

In the present paper I recall the multivariate version of
this result where a stationary random field is considered
with elements indexed by the lattice points $p\in{\mathbb Z}^\nu$.
I do this in an overview about the results in~\cite{10}
and~\cite{11}. In this overview I also speak about the random
spectral measure of a vector valued  Gaussian stationary
random field which yields a spectral representation of the
(vector valued) elements of these random fields. This is a
natural vector valued counterpart of the result about the
spectral representation of scalar valued Gaussian stationary
random fields.

There are also some other notions and results related
to vector valued Gaussian stationary random fields whose
discussion is needed in the proof of the non-central
limit theorem for non-linear functionals of such random
fields. Such notions are the generalized
vector valued Gaussian random fields, their spectral and
random spectral measures, and the multiple Wiener--It\^o
integral with respect to the coordinates of a vector
valued random spectral measure. They are introduced,
and their most important properties are proved in
papers~\cite{10} and~\cite{11}. The goal of this paper
is to give a correct proof of the multivariate generalization
of the results in~\cite{6} with their help.

In short, in my opinion the proof about the multivariate
generalization of the result in~\cite{6} must be
started from the very beginning. First the basic results
about the behavior of vector valued stationary random fields
must be worked out. This is missing from Arcones' paper.
Moreover, it seems to me that a direct proof of Theorem~6
in~\cite{1} would be not more difficult than the proof of
its reduced version presented in formulas~(3.6) and~(3.7)
of~\cite{1}.

Let me remark that although Arcones' proof of the
non-central limit theorem for non-linear functionals
of vector valued Gaussian random fields was problematic,
the proof of its counterpart about the central limit
theorem for such linear functionals under appropriate
conditions was correct. Moreover, in the study of this
result he proved such an estimate in Lemma~1 of
his paper which was applied also in this work.

This paper consists of five sections and two appendices.
In Section~2 the basic notions and results of
papers~\cite{10} and~\cite{11} are recalled. Section~3
contains the main results of this paper. In Section~4 the
preparatory lemmas needed in the proof of the basic theorems
are presented. Section~5 contains the proof of these theorems.
In Appendix~A the background of the limit theorems of this
paper is discussed. In Appendix~B I prove that not only
the finite dimensional distributions of the stochastic
processes considered in Theorem~3.4 converge, but these
processes also weakly  converge to their limit.

\medskip\noindent
{\it Remark:}\/ It was professor Herold Dehling who asked 
me to clarify the proof of Theorem~6 in Arcones' paper~\cite{1}. 
The goal of this work together with the preliminary papers~\cite{10} 
and~\cite{11} was to answer Dehling's question. 
It turned out that to settle this problem first the theory of 
vector valued stationary Gaussian random fields has to be 
worked out. This theory is similar to the theory of  
scalar valued Gaussian random fields, but there are also 
some essential differences between them. Hence the theory of 
vector valued stationary Gaussian random fields cannot be 
considered as a simple generalization of the theory in the 
scalar valued case. I am grateful to professor Dehling for 
calling my attention to this problem.  

\section{On some properties of vector valued Gaussian stationary
random fields.}

In this section I present the most important results
of~\cite{10} and~\cite{11} needed in this paper. At this 
point I do not give their detailed formulation. I shall
present them in a more detailed form when they appear
in our investigation.

We are working with a $d$-dimensional vector valued Gaussian 
stationary random fields $X(p)=(X_1(p),\dots,X_d(p))$, 
$p\in{\mathbb Z}^\nu$, where ${\mathbb Z}^\nu$ denotes 
the lattice points with integer coordinates in the 
$\nu$-dimensional Euclidean space ${\mathbb R}^\nu$ with 
expectation $EX_j(0)=0$ for all $1\le j\le d$. The distribution 
of such random fields is determined by their covariance
function $r_{j,j'}(p)=EX_j(0)X_{j'}(p)=EX_j(m)X_{j'}(m+p)$, 
$1\le j,j'\le d$, $m,p\in{\mathbb Z}^\nu$.

In a result of~\cite{10} it was shown that this covariance
function $r_{j,j'}(p)$, $1\le j,j'\le d$, can be presented in
the following way. For all $1\le j,j'\le d$ there is a complex
measure $G_{j,j'}$ on the torus $[-\pi,\pi)^\nu$ with finite total
variation such that $r_{j,j'}(p)=\int e^{i(p,x)}G_{j,j'}(\,dx)$ for all
$p\in{\mathbb Z}^\nu$, and $G=(G_{j,j'})$, $1\le j,j'\le d$, is an
even, positive semidefinite matrix valued measure on the torus
$[-\pi,\pi)^\nu$. $G$ is called the spectral measure of the random
field $X(p)$, $p\in{\mathbb Z}^\nu$. (A $d$-dimensional matrix
valued measure $G=(G_{j,j'})$, $1\le j,j'\le d$, is called even
if $G_{j,j'}(-A)=\overline{G_{j,j'}(A)}$ for all $1\le j,j'\le d$
and measurable sets~$A$.) For a more detailed discussion see
Section~2 in~\cite{10}.

In Section~3 of~\cite{10} I also defined a $d$-dimensional vector
valued random spectral measure $Z_G=(Z_{G,1},\dots,Z_{G,d})$
corresponding to a $d$-dimensional matrix valued spectral measure
$G$ together with a random integral with respect to it in such a
way that the random integrals $X_j(p)=\int e^{i(p,x)}Z_{G,j}(\,dx)$,
$p\in{\mathbb Z}^\nu$, $1\le j\le d$, define a $d$-dimensional
Gaussian stationary random field with matrix valued spectral
measure~$G$. Besides, I gave the basic properties of a random
spectral measure~$Z_G$ corresponding to a spectral measure~$G$.
These properties determine the distribution of the random
spectral measure as a function of the spectral measure to which
it corresponds.

Once, these results are proved it is not difficult to generalize
them to the case of vector valued Gaussian stationary random
fields defined on the lattice $\frac1K\mathbb Z^\nu$ with some
$K>0$. We define the covariance function
$r_{j,j'}(p)=EX_j((0)X_{j'}(p)=EX_j(m)X_{j'}(p+m)$, $1\le j,j'\le d$,
$p,m\in \frac1K\mathbb Z^\nu$ also in this case. There exists
a spectral measure $G=(G_{j,j'})$, $1\le j,j'\le d$, defined
on the torus $[-K\pi,K\pi)^\nu$ which is a $d$-dimensional
matrix valued even measure, and satisfies the identity
$r_{j,j'}(p)=\int e^{i(p,x)}G_{j,j'}(\,dx)$  for all
$p\in\frac1K\mathbb Z^\nu$ and $1\le j,j'\le d$. There
is also a vector valued random spectral measure
$Z_G=(Z_{G,1},\dots,Z_{G.d})$ corresponding
to this spectral measure $G$ such that
$X_j(p)=\int e^{i(p,x)}Z_{G,j}(\,dx)$, $1\le j\le d$,
$p\in\frac1K\mathbb Z^\nu$, is a vector valued Gausssian
stationary random field on $\frac1K\mathbb Z^\nu$ with
expectation zero and spectral measure~$G$.

It is useful  also to consider vector valued Gaussian stationary
random fields defined in the space $\mathbb R^\nu$. It turned
out that it is even more useful to work with vector valued
generalized Gaussian stationary random fields which can be
considered as their generalization. They are random fields
$$
X(\varphi)=(X_1(\varphi),\dots,X_d(\varphi))
\textrm{ with parameter set }\varphi\in {\cal S}
$$
(instead of ${\mathbb Z}^\nu$ or $\mathbb R^\nu$), where
${\cal S}$ denotes the class of real valued functions in
the $\nu$-dimensional Schwartz space. The definitions
applied in the theory of generalized random fields  were
explained in Section~4 of~\cite{10} together with the
notions needed to understand them.

In paper~\cite{10} generalized vector valued Gaussian
stationary random fields were also constructed, and their
properties were explained. Results similar to those
of~Sections~2 and~3 in~\cite{9} about ordinary vector valued
Gaussian stationary random fields were proved for them.
Generalized vector valued stationary Gaussian random fields
were constructed with the help of their spectral measure which
were also defined.

The spectral measure of a generalized vector valued stationary
random field $X(\varphi)=(X_1(\varphi),\dots,X_d(\varphi))$,  
$\varphi\in {\cal S}$, has properties similar to that of an
ordinary vector valued stationary random field, but there are
some important differences between them. It is a $d\times d$
even, positive definite matrix valued function
$G(A)=(G_{j.j'}(A))$, $1\le j,j'\le d$, $A\subset {\mathbb R}^\nu$,
defined on the bounded, measurable subsets of the $\nu$-dimensional
Euclidean space ${\mathbb R}^\nu$ whose restriction to the measurable
subsets of any finite cube $[-K,K]^\nu$ is a matrix valued measure
with coordinates that are complex measures with finite total variation.
On the other hand, $\sup\limits_A|G_{j,j}(A)|$, where supremum is
taken for all bounded, measurable sets $A$ need not be finite. Only 
the weaker condition 
\begin{equation}
\int (1+|x|)^{-r}G_{j,j}(\,dx)<\infty \;\textrm{ for all }
1\le j\le d \textrm{ with some number } r>0, \label{1.2}
\end{equation}
is imposed. This property is called moderate increase at infinity.

The definition and construction of spectral measures of vector valued
generalized Gaussian stationary random fields was done in the
following way.

Let us consider an even, positive definite matrix valued function
$G(A)=(G_{j.j'}(A))$, $1\le j,j'\le d$, defined on the bounded
and measurable sets $A\subset {\mathbb R}^\nu$ with moderate
increase at infinity, and such that the restriction of its
coordinates to a finite cube $[-K,K]^\nu$ is a complex measure with
finite total variation. If there exists a generalized vector
valued Gaussian stationary random field
$X(\varphi)=(X_1(\varphi),\dots,X_d(\varphi))$, $\varphi\in{\cal S}$,
with the additional property $EX_j(\varphi)=0$, $1\le j\le d$,
for all $\varphi\in{\cal S}$ such that the identity
$$
EX_j(\varphi)X_{j'}(\psi)=\int \tilde \varphi(x)\overline{\tilde\psi}(x)
G_{j.j'}(\,dx), \;\;1\le j,j'\le d, \quad\textrm{for all }
\varphi,\psi\in {\cal S}
$$
holds, where $\,\tilde{}\,$ denotes Fourier transform, and overline
means complex conjugate, then this set of matrix valued functions
$G(A)=(G_{j.j'}(A))$, $1\le j,j'\le d$, is called the spectral
measure of this generalized random field $X(\varphi)$,
$\varphi\in{\cal S}$.

For any set of matrix valued functions $G(A)=(G_{j,j'}(A))$,
$1\le j,j'\le d$, with the above properties there exists a
generalized vector valued Gaussian stationary random field
$X(\varphi)$, $\varphi\in{\cal S}$, with expectation zero
whose covariance function $EX_j(\varphi)X_{j'}(\psi)$,
$\varphi,\psi\in{\cal S}$, satisfies the above conditions.
This means that a matrix valued function with the above
properties is the spectral measure of a generalized, vector
valued Gaussian stationary random field. Moreover, the
distribution of this random field is determined by its
spectral measure.

Given the spectral measure $G=(G_{j,j'}(\cdot))$ of a generalized
random field, such a vector valued random spectral measure
$Z_G=(Z_{G,1},\dots,Z_{G,d})$ can be constructed for which
$X_j(\varphi)=\int\tilde\varphi(x)Z_{G,j}(\,dx)$, $1\le j\le d$,
$\varphi\in{\cal S}$, is a generalized Gaussian stationary random
field with spectral measure $G$ and $EX_j(\varphi)=0$ for all
$1\le j\le d$ and $\varphi\in{\cal S}$. We say that such a random
spectral measure is adapted to the generalized spectral measure~$G$.
The basic properties of the random spectral measures adapted to
a generalized spectral measure also were proved. Their distribution
is determined by the spectral measure to which they are adapted.

The introduction of the random spectral measures corresponding
to the spectral measures of generalized Gaussian stationary
random fields turned out to be useful for us. This class of
random spectral measures is much larger than the class of random
spectral measures corresponding to the spectral measure of
a classical vector valued Gaussian random field. The limit in
the limit theorems of this paper could be expressed by
means of a sum of multiple Wiener--It\^o integrals with
respect to such a random spectral measure.

In the subsequent part of the works~\cite{10} and~\cite{11} my
goal was to give a good representation of those random variables
with finite second moment which are measurable with respect 
to the $\sigma$-algebra generated by the random variables 
of the underlying vector valued random field and to present a
useful formula for their shift transforms. Such results turned
out to be very useful in the study of the limit theorems I was
interested in. A good representation can be given with the help
of multiple Wiener--It\^o integrals with respect to vector valued
random spectral measures introduced in Section~5 of~\cite{10}.

To define multiple Wiener-It\^o integrals I considered the matrix
valued spectral measure $G=(G_{j,j'})$, $1\le j,j'\le d$, of a
$d$-dimensional Gaussian stationary random field, (ordinary or
generalized one), and took a random spectral measure
$Z_G=(Z_{G,1},\dots,Z_{G,d})$ corresponding to it. In Section~5
of~\cite{10} I defined for all $n\ge 1$ and sequences of integers
$j_1,\dots,j_n$ with the property $1\le j_s\le d$ for all $1\le s\le n$
a set $K_{n,j_1,\dots,j_n}=K_{n,j_1,\dots,j_n}(G_{j_1,j_1},\dots,G_{j_n,j_n})$
of complex number valued functions with arguments in ${\mathbb R}^{n\nu}$.
(In the terminology of this paper $K_{n,j_1,\dots,j_n}$ is a subset of
the class of complex valued functions $f(x_1,\dots,x_n)$ of $n$
variables with arguments $x_s\in{\mathbb R}^\nu$, $1\le s\le n$).
I defined the $n$-fold Wiener--It\^o integral
$$
I_n(f|j_1,\dots,j_n)=\int f(x_1,\dots,x_n)
Z_{G,j_1}(\,dx_1)\dots Z_{G,j_n}(\,dx_n)
$$
for the functions $f\in K_{n,j_1,\dots,j_n}$. (The definition of the
set of functions $K_{n,j_1,\dots,j_n}$ is recalled in Section~4 of this
paper before the formulation of Proposition~4A.) Then I proved the
most important properties of these random integrals.

In Section~6 of paper~\cite{10} I proved a technical result, called the
diagram formula about the expression of the product of two multiple
Wiener--It\^o integrals as a sum of multiple Wiener--It\^o integrals. 

These results were exploited in paper~\cite{11}. Here I recalled 
the notion of Wick polynomials which turned out to be a useful tool 
in our investigations. Wick polynomials are natural multivariate
generalizations of Hermite polynomials. Their definition together
with their most important properties was recalled from~\cite{9}
in Section~2 of~\cite{11}. Section~2 of~\cite{11} also contains an
important formula about the expression of Wick polynomials by
means of multiple Wiener--It\^o integrals and another important
formula about the calculation of the shift transforms of a random
variable presented in the form of a sum of multiple Wiener--It\^o
integrals. This made possible to reformulate our limit problems
to limit problems about sums of multiple Wiener--It\^o integrals.
A result in Section~3 of~\cite{11} was proved in order to
investigate such problems. It plays an important role in the
investigation of this paper, hence I recalled it in Proposition~4A
of this paper.

\section{Formulation of the main results.}

In this section I present the main results of this paper.
I shall compare both the formulation of the conditions and
the proof of the results with those appearing in the study
of the analogous results in the scalar valued case. But I
shall refer to~\cite{9} instead of~\cite{6} in this
comparison, because in that work the proofs are worked
out in more detail.

I shall work with such random fields for which
$EX_j(p)=0$ for all $1\le j\le d$ and
$p\in{\mathbb Z}^\nu$. Besides this property I shall
impose two kinds of conditions in this paper.
The first of them deals with the covariance function
$r_{j,j'}(p)=EX_j(0)X_{j'}(p)$, $1\le j,j'\le d$,
$p\in{\mathbb Z}^\nu$, of the vector valued Gaussian
stationary random field $X(p)=(X_1(p),\dots,X_d(p))$,
considered in this paper, the second one with the
function $H(x_1,\dots,x_d)$ which appears in the definition
of the random sums whose limit behavior is investigated.

The following condition is imposed about the covariance
function \hfill\break
$r_{j,j'}(p)=EX_j(0)X_{j'}(p)$.
\begin{equation}
\lim_{T\to\infty}\sup_{p\colon\;p\in{\mathbb Z}^\nu,\,|p|\ge T}
\frac{\left|r_{j,j'}(p)
-a_{j,j'}(\frac p{|p|})|p|^{-\alpha}L(|p|)\right|}
{|p|^{-\alpha}L(|p|)}=0 \label{1.3}
\end{equation}
for all $1\le j,j'\le d$, where $0<\alpha<\nu$, $L(t)$,
$t\ge1$, is slowly varying at infinity, bounded in all finite
intervals, and $a_{j,j'}(t)$ is a real valued continuous
function on the unit sphere
${\cal S}_{\nu-1}=\{x\colon\;x\in {\mathbb R}^\nu,\;|x|=1\}$,
which satisfies the identity $a_{j,j'}(x)=a_{j',j}(-x)$
for all $x\in{\cal S}_{\nu-1}$ and $1\le j,j'\le d$. 

\medskip
I construct a vector valued Gaussian stationary random
field which satisfies relation~(\ref{1.3}). This example
indicates that the covariance functions which satisfy~(\ref{1.3})
have some additional properties, too. These properties will be
discussed in~Appendix~A, because they may help in a better
understanding of the picture about the limit theorems of
this paper.

\medskip\noindent
{\it Example for a stationary random field with covariance
function satisfying relation~(\ref{1.3}).}\/ I shall
construct a stationary random field whose covariance
function satisfies~(\ref{1.3}). I will do this by defining
the spectral measure of such a random field. To do this
I recall some results about the Fourier transform of
generalized functions from the literature.

In the $\nu$-dimensional space ${\mathbb R}^\nu$ the Fourier
transform of the homogeneous function $|x|^\lambda$ (as the Fourier
transform of a generalized function) is $C|p|^{-\lambda-\nu}$
with some coefficient $C=C(\lambda,\nu)>0$. (See the list of
Fourier transforms at the end of the book~\cite{7}.) The value
of this coefficient $C(\lambda,\nu)$ is known, but it
has no importance for us.

On the other hand, if  $u(x)$, $x\in{\mathbb R}^\nu$, is a
sufficiently smooth function, concentrated in a compact
domain, and $u(0)=1$, then the Fourier transform of
$|x|^\lambda u(x)$ equals
$\int e^{i(x,p)} |x|^\lambda u(x)\,dx
=C(\lambda,\nu)|p|^{-\lambda-\nu}(1+o(1))$. 

In the following construction the above property of the
Fourier transform of $|x|^{\lambda}u(x)$ will be exploited.
Define some functions $g_{j.j'}(x)$, $1\le j,j'\le d$,
$x\in [-\pi,\pi)^\nu$, in the following way. Take a
non-negative, smooth function $u(x)$ concentrated in
the cube $[-\pi,\pi]^\nu$ such that $u(-x)=u(x)$, and
$u(0)=1$. Put $g_{j,j}(x)=|x|^{\alpha-\nu}u(x)$, for
$1\le j\le d$, and
$g_{j,j'}(x)=\varepsilon_{j,j'}|x|^{\alpha-\nu}u(x)$ 
for $1\le j,j'\le d$ if $j\neq j'$ with a
sufficiently small real valued coefficient
$\varepsilon_{j,j'}$ such that
$\varepsilon_{j,j'}=\varepsilon_{j',j}$. (One could choose
a complex valued coefficient $\varepsilon_{j,j'}$ too,
but this would demand a more complicated argument.)
I claim that $(g_{j,j'}(x))$, $1\le j,j'\le d$,
$x\in[\-\pi,\pi)^\nu$, with the above defined functions
$g_{j,j'}(\cdot)$ is a spectral density function, and the
covariance function $r_{j,j'}(p)$, $1\le j,j'\le d$,
$p\in{\mathbb Z}^\nu$, of a stationary random field with
this spectral density satisfies relation~(\ref{1.3}) with
$|p|^{-\alpha}$, $L(p)=1$, 
$a_{j,j}\left(\frac p{|p|}\right)=C(\alpha-\nu,\nu)$, and
$a_{j,j'}\left(\frac p{|p|}\right)
=\varepsilon_{j,j'}C(\alpha-\nu,\nu)$ for $j\neq j'$.

Indeed, relation~(\ref{1.3}) holds with such a choice,
because $r_{j,j'}(p)$ is the Fourier transform of
$g_{j,j'}(x)$. We still have to check that $(g_{j,j'}(x))$,
$1\le j,j'\le d$, is a spectral density matrix. The main
point is to show that this matrix is positive definite.
This property holds, since this matrix has the form
$|x|^{\alpha-\nu}u(x)(I+D(\varepsilon))$ with a small
matrix $D(\varepsilon)$, where $I$ denotes the
identity matrix.

Observe that the function
$a_{j,j'}\left(\frac p{|p|}\right)|p|^{-\alpha}L(|p|)
=a_{j,j'}\left(\frac p{|p|}\right)|p|^{-\alpha}$
appearing in formula~(\ref{1.3}) with the functions
defined in the above example is the Fourier transform
of $g^{(0)}_{j,j}(x)=|x|^{\alpha-\lambda}$ if the indices
$j$ and $j'$ of the above function agree, and
$g^{(0)}_{j,j'}(x)=\varepsilon_{j,j'}|x|^{\alpha-\lambda}$ if
$j\neq j'$. Besides, the matrix $(g^{(0)}_{j,j'}(x))$,
$1\le j,j'\le d$, is the spectral density of a vector
valued, stationary, generalized random field. This
spectral density has the homogeneity property
$(g^{(0)}_{j,j'}(tx))=t^{\alpha-\nu}(g^{(0)}_{j,j'}(x))$, $1\le j,j'\le d$,
for all $t>0$. The spectral density $(g_{j.j'}(x))$,
$1\le j,j'\le d$, is in some sense close to this spectral
density $(g^{(0)}_{j,j'}(x))$. In Appendix~A I show that the
spectral measure of a vector valued stationary random field
whose covariance matrix satisfies relation~(\ref{1.3}) has
a similar behavior. It is close in some sense to such a
spectral measure which has some homogeneity property. (This
spectral measure belongs to a generalized random field.)
This homogeneity property has deep consequences in the
theory of the limit theorems we are interested in.

\medskip\noindent
{\it Remark:}\/ There is a natural generalization of the results of
the present paper. One may consider such vector valued stationary
Gaussian random fields, where the partial sums of different
coordinates have a limit with different normalization. They
satisfy limit theorems similar to those of the present paper, but
the different behavior of the different coordinates must be taken
into consideration in the choice of the normalization.

Such more general models were considered in the paper~\cite{14} of
Sanchez de Naranjo, who considered models whose covariance matrices
satisfy a generalized version of relation~(\ref{1.3}). Namely,
they satisfy the relation
$$
r_{j,j'}(p)\sim
|p|^{\alpha_{j,j'}}a_{j,j'}\left(\frac p{|p|}\right)L_{j,j'}(|p|)
$$
with such an exponent $\alpha_{j,j'}$ and slowly varying function
$L_{j,j'}(\cdot)$ which may depend on the indices~$j$ and~$j'$. With
a good choice of these quantities an interesting generalization of
the results of the present paper can be obtained. Such results can
be proved by means of a natural generalization of the arguments of
the present paper, but since this would demand a lot of space and
the introduction of many new quantities I omit their discussion
here.

\medskip
Next I explain the condition imposed on the
function $H(x_1,\dots,x_d)$ that appears in the limit theorems
of this paper. In scalar valued models first the special
case $H(x)=H_k(x)$ was considered, where $H_k(x)$ denotes the 
$k$-th Hermite polynomial with leading coefficient~1. Then
it was shown that our limit problem with a function $H(x)$ 
whose expansion by the Hermite polynomials has the
form $H(x)=\sum_{l=k}^\infty c_l H_l(x)$ with starting index~$k$
in the summation can be simply reduced to the special case 
when $H(x)=c_kH_k(x)$. Similar results will be proved in the 
multivariate case. In this case Wick polynomials take 
the role of the Hermite polynomials. But Wick polynomials
appear in this work only in an implicit way. In the models
studied in this paper the Wick polynomials can be 
simply calculated. Such vector valued Gaussian
stationary random fields 
$X(p)=(X_1(p),\dots,X_d(p))$, $p\in{\mathbb Z}^\nu$, 
are considered whose covariance functions satisfy
besides condition~(\ref{1.3}) also the relation
\begin{eqnarray}
&&EX_j^2(0)=1 \textrm{ for all } 1\le j\le d, \textrm{ and } 
EX_j(0)X_{j'}(0)=0  \nonumber \\
&& \qquad\qquad\qquad \textrm{ if } j\neq j', \;\;1\le j,j'\le d.
\label{1.4}
\end{eqnarray}
First I show that this new condition does not mean a real 
restriction of our problem.

\medskip
Let $X(p)=(X_1(p),\dots,X_d(p))$, $p\in{\mathbb Z}^\nu$,
be a vector valued Gaussian stationary random field 
with expectation $EX_j(p)=0$, $p\in{\mathbb Z}^\nu$, 
$1\le j\le d$, and take the random variables 
$X_1(0),\dots,X_d(0)$ in it. An appropriate number
$1\le d'\le d$ can be chosen, and $d'$ random variables
$X'_j(0)=\sum_{l=1}^dc_{j,l}X_l(0)$, $1\le j\le d'$,
can be defined with appropriate coefficients $c_{j,l}$,
$1\le j\le d'$, $1\le l\le d$, with the following
properties. $EX'_j(0)X'_{j'}(0)=\delta_{j,j'}$, 
$1\le j,j'\le d'$, where $\delta_{j,j'}=0$ if $j\neq j'$, 
and $\delta_{j,j}=1$, and  the random variables 
$X_j(0)$, $1\le j\le d$, can be expressed as the linear
combinations of the random variables $X'_l(0)$,
$1\le l\le d'$, i.e. $X_j(0)=\sum_{l=1}^{d'}d_{j,l}X'_l(0)$
for all $1\le j\le d$ with appropriate coefficients $d_{j,l}$.

Let us define the vector valued random field 
$X'(p)=(X'_1(p),\dots,X'_{d'}(p))$ as
$X'_j(p)=\sum_{l=1}^dc_{j,l}X_l(p)$, $1\le j\le d'$, with
the same coefficients~$c_{j,l}$ as in the definition of
$X'_j(0)$ for all $p\in{\mathbb Z}^\nu$. Then it is not
difficult to see that $X'(p)$, $p\in{\mathbb Z}^\nu$,
is a $d'$-dimensional Gaussian  stationary random field
whose elements have expectation~zero, and it satisfies
relation~(\ref{1.4})  (with parameter~$d'$ instead
of~$d$.) Moreover, if  the covariance function of the
original random field  $X(p)$ satisfied
relation~(\ref{1.3}), then the covariance function of
this new random field also satisfies this condition
with appropriate new functions $a'_{j,j'}(\frac p{|p|})$.
Besides, it is not difficult to find such a function
$H'(x_1,\dots,x_{d'})$ for which
$H'(X'_1(p),\dots,X'_{d'}(p))=H(X_1(p),\dots,X _d(p))$
for all $p\in{\mathbb Z}^\nu$. This means that
with the introduction of this new random field
$X'(p)=(X'_1(p)\dots,X'_{d'}(p))$ our problem can be
reformulated in such a way that our vector valued 
stationary Gaussian random field satisfies both 
relations~(\ref{1.3}) and~(\ref{1.4}). We shall work 
with such a new $d'$-dimensional random field $X'(p)$ 
and function $H'(x_1,\dots,x_{d'})$, only the sign
prime will be omitted everywhere. 

\medskip
First we consider the case when we fix a positive
integer~$k$, and the function $H(x_1,\dots,x_d)$ has the form
\begin{eqnarray}
H(x_1,\dots,x_d)&=&H^{(0)}(x_1,\dots,x_d) \label{1.5} \\
&=&\sum_{\substack{(k_1,\dots, k_d),\; k_j\ge 0,\;1\le j\le d,\\
k_1+\cdots+k_d=k}}
c_{k_1,\dots,k_d}H_{k_1}(x_1)\cdots H_{k_d}(x_d) \nonumber 
\end{eqnarray}
with the previously fixed number~$k$, the coefficients
$c_{k_1,\dots,k_d}$ are real numbers, and $H_{k_j}(\cdot)$ denotes
the Hermite polynomial of order $k_j$ with leading coefficient~1.
The function $H'(x_1,\dots,x_{d'})$ preserves this property of
the function $H(x_1,\dot,x_d)$ when the previously mentioned
transformation is applied, only different coefficients
$c'_{k_1,\dots,k_d}$ appear in its expansion.

\medskip\noindent
{\it Remark:}\/ Although I shall not apply the observation of
this remark, it may be worth mentioning that if $(X_1,\dots,X_d)$
is a $d$-dimensional random vector with standard normal distribution
then $H(X_1,\dots,X_d)$ with a function $H(x_1,\dots,x_d)$ having
the form (\ref{1.5}) is a Wick polynomial of order~$k$ of the
random vector $(X_1,\dots,X_d)$. (See e.g. Corollary~2C in~\cite{11}
or Corollary~2.3 in~\cite{9}.) In general, one can say that Hermite
polynomials play an important role in limit theorems for non-linear
functionals of scalar valued Gaussian random fields. In the case of
non-linear functionals of vector valued Gaussian random fields Wick
polynomials take their role.

\medskip
In scalar valued models, i.e.\ in the case $d=1$ a non-central
limit theorem was proved if $H(x)=H_k(x)$, $k\ge2$, and the 
covariance function $r(n)=EX_0X_n$ satisfies condition~(\ref{1.3})
(with $d=1$) with $0<\alpha<\frac\nu k$. This result was 
formulated in Theorem~8.2 of~\cite{9}. In this result the 
limit was described by means of a $k$-fold Wiener--It\^o
integral with respect to an appropriate random spectral
measure. This random spectral measure corresponds to the
spectral measure that appeared in Lemma~8.1 of~\cite{9} as
the limit of a sequence of appropriately normalized versions
of the spectral measure of a stationary random field $X(p)$,
$p\in{\mathbb Z}^\nu$, whose covariance function satisfies
condition~(\ref{1.3}) with $d=1$. Here I prove a multivariate
version of Theorem~8.2 of~\cite{9} with the help of a multivariate
version of Lemma~8.1 in~\cite{9} formulated below. 

This generalization of Lemma~8.1 in~\cite{9} is a limit
theorem for a sequence of appropriately rescaled versions of
the coordinates $G_{j.j'}$ of a spectral measure $G=(G_{j.j'})$,
$1\le j,j'\le d$, with some nice properties. In this
limit theorem the vague convergence of complex measures is
considered. Before the formulation of this result I recall the
definition of this convergence from Section~3 of~\cite{11}.
In this definition the notion of complex measures with locally
finite total variation appears. I explain its meaning in a
remark after the definition.

\medskip\noindent
{\bf Definition of vague convergence of complex  
measures on ${\mathbb R}^\nu$ with locally finite total variation.} 
{\it Let $G^{(N)}$, $N=1,2,\dots$, be a sequence of complex  
measures on ${\mathbb R}^\nu$ with locally finite total variation. 
We say that this sequence $G^{(N)}$ vaguely converges to a complex
measure~$G^{(0)}$ on~${\mathbb R}^\nu$ with locally finite total 
variation if
$$
\lim_{N\to\infty}\int f(x)\,G^{(N)}(\,dx)=\int f(x)\,G^{(0)}(\,dx) 
$$
for all continuous functions~$f$ on ${\mathbb R}^\nu$ with a 
bounded support.}

\medskip\noindent
{\it Remark:}\/ In the above definition the notion of complex
measures with locally finite total variation appeared. This
notion was introduced in Section~4 of~\cite{10} together with the
notion of vector valued Gaussian stationary generalized random
fields and their matrix valued spectral measures. A complex measure
with locally finite total variation is such a complex valued
function on the bounded measurable subsets of ${\mathbb R}^\nu$
whose restriction to the measurable subsets of a cube $[-T,T]^\nu$
is a complex measure with finite total variation for all $T>0$.

\medskip
The above definition of vague convergence slightly differs from
the classical one presented e.g. in Section~8 of~\cite{9} (before
 Lemma~8.1 of this paper), where the vague convergence of
locally finite (non-negative) measures is considered. The
locally finite measures were defined on all measurable
subsets of ${\mathbb R}^\nu$. Here we deal with complex
measures, because we also want to study the non-diagonal
elements $G_{j,j'}$, $j\neq j'$, of a matrix valued spectral
measure, and they are complex (i.e. not necessary real valued)
measures. A non-negative locally finite measure always can be
extended to a measure on all measurable subsets of
${\mathbb R}^\nu$, while there are locally finite complex
measures which do not have this property. This fact was taken
into account in the introduction of the above definition.

\medskip
The next Proposition~3.1 contains the multivariate version of
Lemma~8.1 in~\cite{9}.

\medskip\noindent
{\bf Proposition~3.1.} {\it Let $G=(G_{j,j'})$ be the matrix
valued spectral measure of a $d$-dimensional vector valued 
stationary random field whose covariance function 
$r_{j,j'}(p)$  satisfies relation~(\ref{1.3}) with some parameter
$0<\alpha<\nu$ and slowly varying function $L(\cdot)$. Let us
define the following rescaled versions of the coordinates
$G_{j,j'}$, $1\le j,j'\le d$, of this matrix valued spectral measure:
\begin{equation}
G^{(N)}_{j,j'}(A)=\frac{N^\alpha}{L(N)}G_{j,j'}\left(\frac AN\right),
\quad A\in{\cal B}^\nu, \quad 1\le j,j'\le d, \label{1.6}
\end{equation}
for all $N=1,2,\dots$,where ${\cal B}^\nu$ denotes the
$\sigma$-algebra of the Borel measurable sets on ${\mathbb R}^\nu$.
Then $G^{(N)}_{j,j'}$ is a complex measure with finite total variation
concentrated in $[-N\pi,N\pi)^\nu$.

For all pairs $1\le j,j'\le d$ the sequence of complex
measures~$G^{(N)}_{j,j'}$ defined in~(\ref{1.6}) tends vaguely to
a complex measure~$G^{(0)}_{j,j'}$ on ${\mathbb R}^\nu$ with locally
finite total variation. These complex measures~$G^{(0)}_{j,j'}$,
$1\le j,j'\le d$, have the homogeneity property
\begin{eqnarray}
&&G^{(0)}_{j,j'}(A)=t^{-\alpha}G^{(0)}_{j,j'}(tA) 
\;\; \textrm{for all bounded sets }A\in{\cal B}^\nu, \nonumber \\ 
&&\qquad\qquad\qquad \;1\le j,j'\le d,\textrm{ and } t>0.  \label{1.7}
\end{eqnarray}

The complex measure $G^{(0)}_{j,j'}$ is determined by the number
$0<\alpha<\nu$ and functions $a_{j,j}(\cdot)$, $a_{j,j'}(\cdot)$,
$a_{j',j}(\cdot)$ and $a_{j',j'}(\cdot)$ defined in
formula~(\ref{1.3}) on the unit sphere~$S_{\nu-1}$. This
implies that for all spectral measures $G$ that satisfy
relation~(\ref{1.3}) with the same parameter $\alpha$ and
functions $a_{j,j'}(\cdot)$, $1\le j,j'\le d$ the vague limit
of the complex measures $G_{j.j'}^{(N)}$ is the same for all
$1\le j,j'\le d$.

Finally, there exists a vector valued generalized Gaussian
stationary random field on ${\mathbb R}^\nu$ whose matrix valued
spectral measure is  $G^{(0)}=(G^{(0)}_{j,j'})$, $1\le j,j'\le d$,
with the complex measures $G^{(0)}_{j,j'}$, $1\le j,j'\le d$,
defined in this Proposition.} 

\medskip
In the following Theorem~3.2 I formulate the multivariate
version of Theorem~8.2 in~\cite{9}. In its formulation the
result of Proposition~3.1 is applied where a matrix valued
spectral measure $G^{(0)}=(G^{(0)}_{j.j'})$, $1\le j,j'\le d$, is
constructed under some conditions which are imposed also in
Theorem~3.2. Theorem~3.2 is a limit theorem where the limit
is defined by means of a sum of multiple Wiener--It\^o integrals
with respect to a vector valued random spectral measure that
corresponds to the matrix valued spectral measure $G^{(0)}$
constructed in Proposition~3.1. Let me remark that I formulated
this result also in paper~\cite{10}. But in that work it was not
proved. That work contained only a heuristic argument which
indicated why it is natural to expect such a result. Its goal
was to indicate the usefulness of the theory worked out
in~\cite{10} and \cite{11}.

\medskip\noindent
{\bf Theorem~3.2.} {\it Fix an integer $k\ge1$,
and let $X(p)=(X_1(p),\dots,X_d(p))$, 
$p\in{\mathbb Z}^\nu$, be a vector valued Gaussian 
stationary random field whose covariance matrix 
$r_{j,j'}(p)=EX_j(0)X_{j'}(p)$, $1\le j,j'\le d$, 
$p\in{\mathbb Z}^\nu$, satisfies both relation~(\ref{1.3})
with some number $\alpha$ such that $0<\alpha<\frac\nu k$
and relation~(\ref{1.4}). Let  $H(x_1,\dots,x_d)$ be a
function of the form given in~(\ref{1.5}) also with the
previously fixed number~$k$. Define the random variables 
$Y(p)=H(X_1(p),\dots,X_d(p))$ for all 
$p\in{\mathbb Z}^\nu$ together with their normalized 
partial sums 
$$
S_N=\frac1{N^{\nu-k\alpha/2}L(N)^{k/2}}\sum_{p\in B_N}Y(p),
$$ 
where the set $B_N$ was defined in~(\ref{1.1}). These random 
variables $S_N$, $N=1,2,\dots$, satisfy the following 
limit theorem.

Let $Z_{G^{(0)}}=(Z_{G^{(0)},1},\dots,Z_{G^{(0)},d})$ be a vector 
valued random spectral measure which corresponds to the 
matrix valued spectral measure $G^{(0)}=(G^{(0)}_{j,j'})$, 
$1\le j,j'\le d$, defined in Proposition~3.1 with
the help of the matrix valued spectral measure $G=(G_{j,j'})$
of a vector valued Gaussian stationary random field with
covariance function $r_{j,j'}(s)$, $1\le j,j'\le d$,
$s\in{\mathbb R}^\nu$, satisfying relation~(\ref{1.3}).
Then the sum of multiple Wiener--It\^o integrals with
the coefficients $c_{k_1,\dots,k_d}$ appearing in~(\ref{1.5})
\begin{eqnarray}
S_0&=&\sum_{\substack{(k_1,\dots, k_d),\; k_j\ge 0,\;1\le j\le d, \\
k_1+\cdots+k_d=k}} c_{k_1,\dots,k_d}
\int \prod_{l=1}^\nu\frac{e^{i(x_1^{(l)}+\cdots+x_k^{(l)})}-1}
{i(x^{(l)}_1+\cdots+x_k^{(l)})} \label{1.8} \\
&&\qquad\qquad\qquad Z_{G^{(0)},{j(1|k_1,\dots,k_d)}}(\,dx_1)\dots 
Z_{G^{(0)},{j(k|k_1,\dots,k_d)}}(\,dx_k)
\nonumber
\end{eqnarray}
exists, where the notation  
$x_p=(x_p^{(1)},\dots,x_p^{(\nu)})\in{\mathbb R}^\nu$,
$p=1,\dots,k$, is applied, and the indices
$j(s|k_1,\dots,k_d)$, $1\le s\le k$, are defined as
$j(s|k_1,\dots,k_d)=r$ if
$\sum_{u=1}^{r-1} k_u< s\le \sum_{u=1}^r k_u$, $1\le s\le k$,
$1\le r\le d$. (For $r-1=0$ the convention
$\sum_{u=1}^0 k_u=0$ is applied in this definition.)
The normalized sums~$S_N$ converge in distribution to
the random variable $S_0$ defined in~(\ref{1.8})
as $N\to\infty$.}

\medskip
The indexation of the terms $Z_{G^{(0)},{j(s|k_1,\dots,k_d)}}(\,dx_s)$
in formula (\ref{1.8}) can be described in a simpler form. In
the first $k_1$ arguments $x_1,\dots,x_{k_1}$, i.e. for
$1\le s\le k_1$ $Z_{G^{(0)},1}(\,dx_s)$, is written, in
the next $k_2$ arguments, i.e. for $k_1+1\le s\le k_1+k_2$ 
$Z_{G^{(0)},2}(\,dx_s)$ is written, and so on. In the last $k_d$
arguments, i.e, when $k_1+\cdots+k_{d-1}+1\le s\le k$,
($k=k_1+\cdots+k_d$), $Z_{G^{(0)},d}(\,dx_s)$ is written.

\medskip
In Theorem~3.2 the limit of 
$A_N^{-1}\sum_{p\in B_N} H(X_1(p),\dots,X_d(p))$ is described if
the expansion of the function $H(x_1,\dots,x_d)$ is a linear
combination of products of Hermite polynomials with
different arguments, and all these products are polynomials
of order~$k$. The next Theorem~3.3 which is the multivariate
version of Theorem~8.2$'$ in~\cite{9} states that a similar
result holds if the function $H(x_1,\dots,x_d)$ is the linear
combination of products of Hermite polynomials, but some of
these products may be polynomials of order higher than $k$.

\medskip\noindent
{\bf Theorem~3.3.} {\it Let us consider a vector 
valued Gaussian stationary random field 
$X(p)=(X_1(p),\dots,X_d(p))$, $p\in{\mathbb Z}^\nu$, that
satisfies the conditions of Theorem~3.2 and a function of the form 
$H(x_1,\dots,x_d)=H^{(0)}(x_1,\dots,x_d)+H^{(1)}(x_1,\dots,x_d)$,
where $H^{(0)}(x_1,\dots,x_d)$ was defined in~(\ref{1.5}), and
\begin{equation}
H^{(1)}(x_1,\dots,x_d)
=\sum_{\substack{(k_1,\dots, k_d),\; k_j\ge 0,\;1\le j\le d,\\
k_1+\cdots+k_d\ge k+1}} 
c_{k_1,\dots,k_d}H_{k_1}(x_1)\cdots H_{k_d}(x_d) \label{1.9}
\end{equation}
with real valued coefficients $c_{k_1,\dots,k_d}$ such that
\begin{equation}
\sum_{\substack{(k_1,\dots, k_d),\; k_j\ge 0,\;1\le j\le d,\\
k_1+\cdots+k_d\ge k+1}} 
\frac{c^2_{k_1,\dots,k_d}}{k_1!\cdots k_d!}<\infty. \label{1.10}
\end{equation}
Define the random variables 
$Y(p)=H(X_1(p),\dots,X_d(p))$ for all $p\in{\mathbb Z}^\nu$
and their normalized partial sums 
$$
S_N=\frac1{N^{\nu-k\alpha/2}L(N)^{k/2}}\sum_{p\in B_N}Y(p),\quad
N=1,2,\dots,
$$ 
with this function $H(x_1,\dots,x_d)$. The 
random variables $S_N$ converge in distribution to the 
random variable $S_0$ defined in formula~(\ref{1.8}) as $N\to\infty$.}

\medskip
Actually condition~(\ref{1.10}) in Theorem~3.3 means that 
$$
E\left[{H^{(1)}}(X_1(0),\dots,X_d(0))^2\right]<\infty.
$$ 
Finally I mention that Arcones formulated a more general 
result. To present it, more precisely to present its
generalization to the case when we are working with
stationary random fields parametrized by the lattice
points of ${\mathbb Z}^\nu$ with some $\nu\ge1$ let us define
the following parameter sets for all $N=1,2,\dots$ and 
$t=(t_1,\dots,t_\nu)$, $0\le t_l\le 1$, for all $1\le l\le \nu$.
\begin{eqnarray}
 B_N(t)&=&B_N(t_1,\dots,t_\nu) \label{1.11} \\
&=&\{p=(p_1,\dots,p_\nu)\in{\mathbb Z}^\nu\colon\; 0< p_l\le Nt_l
\textrm { for all } 1\le l\le\nu\}. \nonumber
\end{eqnarray}
With this notation the following result holds.

\medskip\noindent
{\bf Theorem~3.4.} {\it Let us consider the same vector 
valued Gaussian stationary random field 
$X(p)=(X_1(p),\dots,X_d(p))$, $p\in{\mathbb Z}^\nu$, 
and function $H(x_1,\dots,x_d)$ as in Theorem~3.3.
Define the random variables 
$Y(p)=H(X_1(p),\dots,X_d(p))$ for all $p\in{\mathbb Z}^\nu$
together with the random fields 
\begin{equation}
S_N(t)=\frac1{N^{\nu-k\alpha/2}L(N)^{k/2}}\sum_{p\in B_N(t)}Y(p)
\label{1.12}
\end{equation}
with parameter set $t=(t_1,\dots,t_\nu)$, $0\le t_l\le 1$,
$1\le l\le\nu$, for all $N=1,2,\dots$, where the set $B_N(t)$
was defined in~(\ref{1.11}). The finite dimensional 
distributions of the random fields $S_N(t)$ converge to 
that of the random field $S_0(t)$, $t=(t_1,\dots,t_\nu)$, 
$0\le t_l\le 1$, $1\le l\le\nu$, defined by the formula
\begin{eqnarray}
S_0(t)&=&\sum_{\substack{(k_1,\dots, k_d),\; k_j\ge 0,\;1\le j\le d,\\
k_1+\cdots+k_d=k}}
c_{k_1,\dots,k_d}
\int \prod_{l=1}^\nu\frac{e^{it_l(x_1^{(l)}+\cdots+x_k^{(l)})}-1}
{i(x^{(l)}_1+\cdots+x_k^{(l)})} \label{1.13} \\
&&\qquad\qquad Z_{G^{(0)},{j(1|k_1,\dots,k_d)}}(\,dx_1)\dots 
Z_{G^{(0)},{j(k|k_1,\dots,k_d)}}(\,dx_k)  \nonumber
\end{eqnarray}
if the limit $N\to\infty$ is taken. Similarly to Theorem~3.2 
the notation $x_p=(x_p^{(1)},\dots,x_p^{(\nu)})$, $p=1,\dots,k$,
is applied, and  the indices $j(s|k_1,\dots,k_d)$,
$1\le s\le k$, are defined in the same way as in
formula~(\ref{1.8}).}

\medskip
A referee proposed to show that also a strengthened form
of Theorem~3.4 formulated in the next Corollary holds. I
shall present the proof of this result in Appendix~B. I shall
omit some technical details of the proof, and in the case
$\nu>1$ I shall apply a result whose formulation I did not
find in the literature. I chose such an approach, because a
detailed proof would demand the elaboration of many
complicated technical details which are not related to the
subject of this paper.

\medskip\noindent
{\bf Corollary of Theorem~3.4.} {\it Under the conditions of
Theorem~3.4 not only the finite dimensional distributions of
the random fields $S_N(t)$, $t=(t_1,\dots,t_\nu)$,
$0\le t_l\le1$, $1\le l\le \nu$, introduced in~(\ref{1.12})
converge to those of the random field $S_0(t)$ defined
in~(\ref{1.13}), but even the distribution of the random fields
$S_N(\cdot)$, converge weakly to the distribution of $S_0(\cdot)$
in the Skorochod space on $[0,1]^\nu$ as $N\to\infty$. Moreover,
the trajectories of $S_0(\cdot)$ are continuous functions
on $[0,1]^\nu$.}

\medskip
Let us observe that the kernel functions in the Wiener--It\^o
integrals appearing in the sum which defines $S_0(t)$ in~(\ref{1.13})
equal $\varphi_t(x_1+\cdots+x_k)$, where $\varphi_t(u)$,
$u\in {\mathbb R}^\nu$, is the Fourier transform of the Lebesgue
measure on the rectangle $[0,t_1]\times\cdots\times[0 ,t_\nu]$. The
integral in~(\ref{1.13}) is taken on the whole space.

\medskip
Theorem~3.4 was formulated in that form as Arcones did, but it
could have been formulated in a slightly more general form.
The sets $B_{N}(t)$ in~(\ref{1.11}), the random variables $S_N(t)$,
$N=1,2,\dots$, in~(\ref{1.12}) and $S_0(t)$ in~(\ref{1.13}) could
have been defined for all $t=(t_1,\dots,t_\nu)\in[0,\infty)^\nu$
and not only for $t=(t_1,\dots,t_\nu)\in[0,1]^\nu$. After the
introduction of these objects it could have been proved, similarly to
the proof of Theorem~3.4, that the finite dimensional distributions
of the random fields $S_N(t)$ converge to the finite dimensional
distributions of the random field $S_0(t)$ as $N\to\infty$ also
in this more general case. This more general form of the result is
useful, because it makes possible to formulate an important
property of the limit field $S_0(t)$, called the self-similarity
property. The limit random field $S_0(t)$, $t\in[0,\infty)^\nu$,
is self-similar with parameter $\nu-k\alpha/2$, which means that
$S_0(ut)\stackrel{\Delta}{=}u^{\nu-k\alpha/2}S_0(t)$ for all $u>0$,
where $\stackrel{\Delta}{=}$ means that the finite dimensional
distributions of the two random fields agree.

The self-similarity property of the random field $S_0(t)$,
$t\in[0,\infty)^\nu$, can be proved by exploiting that
by formula~(\ref{1.7}) in Proposition~3.1 $G^{(0)}(uA)=u^\alpha G^{(0)}(A)$
for the spectral measure $G^{(0)}$ for all $u>0$ and measurable
sets $A\subset{\mathbb R}^\nu$. This implies that 
$$
(Z_{G^{(0)},1}(uA_1),\dots,Z_{G^{(0)},d}(uA_d))\stackrel{\Delta}{=}
(u^{\alpha/2}Z_{G^{(0)},1}(A_1),\dots,u^{\alpha/2}Z_{G^{(0)},d}(A_d))
$$
for all $u>0$ and measurable sets 
$A_1\in{\mathbb R}^\nu,\dots,A_d\in{\mathbb R}^\nu$.
We still have to exploit that the kernel functions
$$
f_t(x_1,\dots,x_k)=\prod_{l=1}^\nu\frac{e^{it_l(x_1^{(l)}+\cdots+x_k^{(l)})}-1}
{i(x^{(l)}_1+\cdots+x_k^{(l)})}
$$ 
in the Wiener--It\^o integrals in~(1.13) (with the notation 
$t=(t_1,\dots,t_\nu)$) have the property 
$$
f_{ut}(x_1,\dots,x_k)=u^\nu f_t(ux_1,\dots,ux_k)
$$ 
for all $u>0$, $t\in[0,\infty)^\nu$, $x_j\in{\mathbb R}^\nu$, 
$1\le j\le d$. The self-similarity property of the random
field $S_0(t)$, $t\in[0,\infty)^\nu$, can be proved with
the help of the above observations.

\section{Preparatory results for the proof of the main theorems.}

This section contains the proof of Proposition~3.1 and the
elaboration of a method that helps in proving the theorems
of this paper. In the application of this method the normalized
random sums $S_N$ appearing in the formulation of Theorem~3.2
are rewritten in the form of a sum of multiple Wiener-It\^o
integrals with respect to a vector valued random spectral
measure. Then Proposition~3.1 of paper~\cite{11} is recalled,
and it is shown how the sums of Wiener--It\^o integrals
expressing the random sums~$S_N$ can be investigated with
its help.

First I prove Proposition~3.1.

\medskip\noindent
{\it Proof of Proposition~3.1.} Proposition~3.1 is proved by means
of an adaptation of the proof of Lemma~8.1 in~\cite{9}. The same
argument works, only some steps of the proof must be modified in a
natural way. I do not work out all details, I only briefly remark what
kind of modifications are needed.

The  diagonal elements $G_{j,j}$, $1\le j\le d$, of the matrix 
valued spectral measure $G$ are spectral measures. Hence
Lemma~8.1 of~\cite{9} implies that for any $1\le j\le d$ the
measures $G^{(N)}_{j,j}$ converge vaguely to a locally finite
measure $G^{(0)}_{j,j}$ determined by the function $a_{j,j}(\cdot)$ and
the number~$\alpha$ which appears in relation~(\ref{1.7}).

For the non-diagonal elements $G_{j,j'}$, $j\neq j'$, this argument
cannot be applied, because $G_{j,j'}$ is a complex measure with
finite total variation which may be not a (positive) measure.
In this case it can be exploited that $G$ is a positive semidefinite
matrix valued measure. Hence the $2\times2$ matrix
$$
G(A|j,j')=\left( \begin{array}{l}
G_{j,j}(A),\quad\,\, G_{j,j'}(A) \\
G_{j',j}(A),\quad G_{j',j'}(A)
\end{array} \right)
$$
is positive semidefinite for all pairs $1\le j,j'\le d$, $j\neq j'$,
and measurable sets $A\subset {\mathbb R}^\nu$. This implies that the 
quadratic forms
$$
(1,1)G(A|j,j')(1,1)^*=G_{j,j}(A)+G_{j',j'}(A)+G_{j,j'}(A)+G_{j',j}(A)
$$
and
$$
(1,i)G(A|j,j')(1,-i)^*=G_{j,j}(A)+G_{j',j'}(A)-i[G_{j,j'}(A)-G_{j',j}(A)]
$$
are non-negative numbers for all measurable sets 
$A\subset {\mathbb R}^{\nu}$. Therefore the set-functions 
$R_{j,j'}(\cdot)$ and $S_{j,j'}(\cdot)$ defined as 
$R_{j,j'}(A)=G_{j,j}(A)+G_{j',j'}(A)+G_{j,j'}(A)+G_{j',j}(A)$ and
$S_{j,j'}(A)=G_{j,j}(A)+G_{j',j'}(A)-i[G_{j,j'}(A)-G_{j',j}(A)]$ for 
all measurable sets~$A\in[-\pi,\pi)^\nu$ are finite measures. 
Their Fourier transforms equal  
$r^{(1)}_{j,j'}(p)=\int e^{i(p,x)}R_{j,j'}(\,dx)=r_{j,j}(p)+r_{j',j'}(p)
+r_{j.j'}(p)+r_{j',j}(p)$ and
$r^{(2)}_{j,j'}(p)=\int e^{i(p,x)}S_{j,j'}(\,dx)=r_{j,j}(p)+r_{j',j'}(p)
+i[r_{j.j'}(p)-r_{j',j}(p)]$, $p\in\mathbb Z^\nu$. These Fourier
transforms satisfy the following relation, similar to formula~(\ref{1.3}).
\begin{equation}
\lim_{T\to\infty}\sup_{p\colon\;p\in{\mathbb Z}^\nu,\,|p|\ge T}
\frac{\left|r^{(s)}_{j,j'}(p)-a^{(s)}_{j,j'}
(\frac p{|p|})|p|^{-\alpha}L(|p|)\right|}
{|p|^{-\alpha}L(|p|)}=0 \label{2.6}
\end{equation}
both for $s=1$ and $s=2$ with some functions $a_{j,j'}^{(s)}(\cdot)$
which can be expressed by means of the functions $a_{j,j}(\cdot)$,
$a_{j.j'}(\cdot)$, $a_{j',j}(\cdot)$ and $a_{j',j'}(\cdot)$. The only
difference from formula~(\ref{1.3}) is that the continuous
function, $a^{(2)}_{j,j'}(\cdot)$ may be complex valued. (I also remark
that the symmetry property $a_{j,j'}(u)=a_{j',j}(-u)$ yields that
$a^{(2)}_{j,j'}(-u)=\overline{a^{(2)}_{j',j}(u)}$. On the other hand,
$a^{(1)}_{j,j'}(\cdot)$ is a real valued function, for which
$a^{(1)}_{j,j'}(-u)=a^{(1)}_{j,j}(u)$. These relations correspond to
the fact that $r^{(s)}_{j,j'}(p)$, $s=1,2$ are Fourier transforms
of real valued measures.)

A natural adaptation of the proof of Lemma~8.1 in~\cite{9} shows
that the measures $R_{j.j'}(\cdot)$ and $S_{j,j'}(\cdot)$ have properties
similar to $G_{j,j}(\cdot)$, only the function $a_{j,j}(\cdot)$ must
be replaced by $a^{(1)}_{j,j'}(\cdot)$ and $a^{(2)}_{j,j'}(\cdot)$ in
them. Moreover, the proof of Lemma~8.1 in~\cite{9} can be applied
to show this. To understand this let us remark that $R_{j,j'}(\cdot)$
and $S_{j,j'}(\cdot)$ are measures on the torus $[-\pi.\pi)^\nu$, and
their Fourier transforms satisfy relation~(\ref{2.6}). The
spectral measure $G(\cdot)$ investigated in Lemma~8.1 of~\cite{9}
has similar properties, and the proof was based on them.

More explicitly, define the measures $R^{(N)}_{j,j'}(\cdot)$,
and $S^{(N)}_{j,j'}(\cdot)$ as
$$
R^{(N)}_{j,j'}(A)=\frac{N^\alpha}{L(N)}R_{j,j'}\left(\frac AN\right), \qquad
S^{(N)}_{j,j'}(A)=\frac{N^\alpha}{L(N)}S_{j,j'}\left(\frac AN\right)
$$
for all measurable sets $A\subset [-N\pi,N\pi)^\nu$ and
$N=1,2,\dots$. I claim that these measures converge vaguely
to some locally finite measures $R^{(0)}_{j.j'}(\cdot)$ and
$S^{(0)}_{j,j'}(\cdot)$ with some homogeneity property
on $\mathbb R^\nu$.

To prove this homogeneity property let us introduce, similarly
to the proof of Lemma~8.1 in~\cite{9} the measures
$\mu^{(1)}_N$ and $\mu^{(2)}_N$, $N=1,2,\dots$ as
$$
\mu^{(1)}_N(A)=\int_A |K_N(x)|^2\,R^{(N)}_{j,j'}(\,dx),
\quad A\in {\cal B}^\nu, \quad N=1,2,\dots
$$
and
$$
\mu^{(2)}_N(A)=\int_A |K_N(x)|^2\,S^{(N)}_{j,j'}(\,dx),
\quad A\in {\cal B}^\nu, \quad N=1,2,\dots
$$
with the function $K_N(\cdot)$ defined as
$$
K_N(x)=\frac1N\sum_{p\in\mathbb B_N}e^{i(p,x/N)}=
\prod_{j=1}^\nu\frac {e^{ix^{(j)}}-1}{N(e^{ix^{(j)}/N}-1)},
\quad N=1,2,\dots.
$$
I claim that both for $s=1$ and $s=2$ the sequence of measures
$\mu^{(s)}_N$ converge weakly to a measure $\mu^{(s)}_0$ as
$N\to\infty$ whose Fourier transform depends on the function
$a^{(s)}_{j,j'}(\cdot)$ and parameter $\alpha$ appearing in
formula~(\ref{2.6}).

This can be proved by calculating the Fourier transforms of the
measures $\mu^{(s)}_N$ and by showing that they have a limit which
is a continuous function. More precisely, Lemma~8.4 in~\cite{9}
must be applied, which yields a modified version of this method. The
reason for this modification is that we can calculate the Fourier
transforms of the measures $\mu^{(s)}_N$, $s=1,2$, only in the points
$\frac pN$, $p\in{\mathbb Z}^\nu$. On the other hand,
the measures $\mu^{(s)}_N$ are concentrated in the cube
$[-N\pi,N\pi)^\nu$. Lemma~8.4 in~\cite{9} provides such a version
of the characteristic function method which works in such cases.

The calculations needed to prove the above properties of the
measures $\mu^{(s)}_N$, $s=1,2$, are carried out in the proof of
Theorem~8.2 in~\cite{9}. Actually, a more general result
is proved there. I omit the details.

Let us define the continuous function
$K_0(x)=\prod_{j=1}^\nu\frac {e^{ix^{(j)}}-1}{ix^{(j)}}$ on $\mathbb R^\nu$.
In all compact subsets of $\mathbb R^\nu$ the functions $K_N(x)$
converge to $K_0(x)$ in the supremum norm as $N\to\infty$. The
proof of Lemma~8.1 in~\cite{9} shows on the basis of this
property that the limit measures $\mu^{(s)}_0$, $s=1,2$, have the
following representation. There are measures $H^{(0)}_{j,j'}$
and $K^{(0)}_{j,j'}$ such that
$\mu^{(1)}_0(A)=\int_A |K_0(x)|^2 H^{(0)}_{j,j'}(\,dx)$, and
$\mu^{(2)}_0(A)=\int_A |K_0(x)|^2 K^{(0)}_{j,j'}(\,dx)$ for all
measurable sets $A\subset \mathbb R^\nu$. Moreover, $H^{(0)}_{j,j'}$ 
and $K^{(0)}_{j,j'}$ are locally finite measures, and they are
the vague limits of the sequences of measures $R^{(N)}_{j,j'}$ and
$S^{(N)}_{j,j'}$ respectively. The measures
$H^{(0)}_{j,j'}$ and $K^{(0)}_{j,j'}$
are determined by the limit measures $\mu^{(1)}_0$ and $\mu^{(2)}_0$,
hence also by the parameter~$\alpha$ and functions
$a_{j,j'}^{(1)}(\cdot)$ and $a_{j,j'}^{(2)}(\cdot)$ in(\ref{2.6}).
The argument of the proof in Lemma~8.1 of~\cite{9} enables us
to show that (\ref{1.7}) holds if the complex measure $G^{(0)}_{j,j}$
is replaced by the measure $H^{(0)}_{j,j'}$ or $K^{(0)}_{j,j'}$ in it.

Since the complex measure $G_{j,j'}$ can be expressed as a linear
combination of the measures $G_{j,j}$, $G_{j',j'}$, $R_{j,j'}$  and
$S_{j,j'}$ the properties proved for them imply the statements
formulated about the behavior of~$G_{j,j'}$ in Proposition~3.1.

We still have to show that $(G^{(0)}_{j,j'})$, $1\le j,j'\le d$,
is the spectral measure of a generalized vector valued stationary
Gaussian random field. By Theorem~4.1 of~\cite{10}  
$(G^{(0)}_{j,j'})$, $1\le j,j'\le d$, is the spectral measure of a
vector valued generalized Gaussian stationary random field if
it is a positive definite matrix valued even measure on
$\mathbb R^\nu$ whose distribution is moderately increasing
at infinity, i.e it satisfies relation~(\ref{1.2}). It follows
from Lemma~3.2 in~\cite{11} and the already proved part of
Proposition~3.1 that this system is a positive semidefinite matrix
valued even measure on ${\mathbb R}^\nu$. The validity of 
relation~(\ref{1.2}) follows from the fact that $G^{(0)}_{j,j}$
has locally finite total variation, and it satisfies
relation~(\ref{1.7}). Proposition~3.1 is proved. 

\medskip
Now I turn to the representation of the normalized random
sums in the form of a sum of multiple Wiener--it\^o integrals.
To do this let us first consider the random variable
$Y(0)=H(X_1(0),\dots,X_d(0))$ defined with the help of
the function $H(x_1,\dots,x_d)=H^{(0)}(x_1,\dots,x_d)$
introduced in~(\ref{1.5}) and a vector valued Gaussian
stationary random field $X(p)=(X_1(p),\dots,X_d(p))$,
$p\in{\mathbb Z}^\nu$, with covariance function 
$r_{j,j'}(p)=EX_j(0)X_{j'}(p)$, $1\le j,j'\le d$, 
$p\in{\mathbb Z}^\nu$, which satisfies relation~(\ref{1.3})
with some parameter $0<\alpha<\frac\nu k$ together with
its shifts $Y(p)=H(X_1(p),\dots,X_d(p))$, $p\in{\mathbb Z}^\nu$,
and express them as a sum of Wiener--It\^o integrals.

This will be done with the help of the results
in~\cite{10} and~\cite{11}.

Let $G=(G_{j,j'})$, $1\le j,j'\le d$, be the matrix valued
spectral measure of the stationary random field
$X(p)=(X_1(p),\dots,X_d(p))$, $p\in{\mathbb Z}^\nu$,  
and let us consider that vector valued random spectral
measure $Z_G=(Z_{G,1},\dots,Z_{G,d})$ corresponding to this
spectral measure for which $X_j(p)=\int e^{i(p,x)}Z_{G,j}(\,dx)$
for all $p\in{\mathbb Z}^\nu$ and $1\le j\le d$. By the
results of~\cite{10} there exists such a vector valued random
spectral measure.

The random variable $Y(0)=H(X_1(0),\dots,X_d(0))$ will be
rewritten in the form of a sum of Wiener--It\^o integrals
with the help of the multiple version of It\^o's formula
presented in Theorem~2.2 of~\cite{11}, more precisely
by the corollary of this result. As $Y(0)$ is a Wick
polynomial of the (independent) random variables $X_j(0)$
with standard Gaussian distribution, and
$X_j(0)=\int Z_{G,j}(\,dy)$, $1\le j\le d$, this formula
yields the desired expression for $Y(0)$. Let me remark
that by Lemma~8B of~\cite{9} condition~(\ref{1.3})
implies that the diagonal elements $G_{j,j}$, $1\le j\le d$,
of the matrix valued spectral measure $G=(G_{j,j'})$,
$1\le j,j'\le d$, are non-atomic. Hence the multiple
Wiener--It\^o integrals with respect to the
coordinates of the vector valued random spectral
measure $Z_G=(Z_{G,1},\dots,Z_{G,d})$ whose sum
expresses $Y(0)$ in the next formula are meaningful.

The above results yield the identity
\begin{eqnarray*}
Y(0)&=&H(X_1(0),\dots,X_d(0))=H^{(0)}(X_1(0),\dots,X_d(0)) \\
&=&\sum_{\substack{(k_1,\dots, k_d),\; k_j\ge 0,\;1\le j\le d,\\
k_1+\cdots+k_d=k}} 
\colon c_{k_1,\dots,k_d} X_1(0)^{k_1}\cdots X_d(0)^{k_d} \colon \\
&=&\sum_{\substack{(k_1,\dots,k_d),\; k_j\ge0\;1\le j\le d,\\ 
k_1+\cdots+k_d=k}} c_{k_1,\dots,k_d} \int \prod_{j=1}^d
\left(\prod_{s=k_1+\cdots+k_{j-1}+1}^{k_1+\cdots+k_j} Z_{G,j}(\,dy_s)\right),
\end{eqnarray*}
where for $j=1$ we define
$\prod\limits_{s=k_1+\cdots+k_{j-1}+1}^{k_1+\cdots+k_j} Z_{G,j}(\,dy_s)
=\prod\limits_{s=1}^{k_1} Z_{G,1}(\,dy_s)$, and if $k_j=0$ for some 
$1\le j\le d$, then we drop the term 
$\prod\limits_{s=k_1+\cdots+k_{j-1}+1}^{k_1+\cdots+k_j}Z_{G,j}(\,dy_s)$
from this expression. (Here 
$\colon\! P(X_1(0),\dots,X_d(0))\colon$ denotes the
Wick polynomial corresponding to $P(X_1(0),\dots,X_d(0))$, 
where $P(x_1,\dots,x_d)$ is a homogeneous polynomial.)

Since $Y(p)=T_pY(0)$ with the shift transformation $T_p$ for all 
$p\in{\mathbb Z}^\nu$, the previous
identity and Proposition~2.4 in~\cite{11} yield the formula
\begin{eqnarray*}
Y(p)&=&T_pY(0) \\
&=& \!\!\!\!\!\!\!\!\!\!\!\!\!\!
\sum_{\substack{(k_1,\dots,k_d),\;k_j\ge0,\;1\le j\le d,\\ 
k_1+\cdots+k_d=k}}  \!\!\!\!\!\!\!\!\!\!\!\!\!
 c_{k_1,\dots,k_d}\int e^{i(p,y_1+\cdots+y_k)}\prod_{j=1}^d
\left(\prod_{s=k_1+\cdots+k_{j-1}+1}^{k_1+\cdots+k_j} \!\!\! 
Z_{G,j}(\,dy_s)\right)
\end{eqnarray*}
for all $p\in{\mathbb Z}^\nu$. By summing up 
this formula for all $p\in B_N$ we get that
\begin{eqnarray*}
&&\!\!\!\!\!\!\!\!\!\!\!\!\!\!\!\!\!\!\!\!\!\!\!\!\!\!\!\!\!
S_N=\frac1{N^{\nu-k\alpha/2}L(N)^{k/2}}
\sum_{\substack{(k_1,\dots,k_d),\;k_j\ge0,\;1\le j\le d\\ k_1+\cdots +k_d=k}}
 c_{k_1,\dots,k_d} \\
&& \!\!\!\!\!\!\!\!\!\!
\int\prod _{l=1}^\nu \frac{e^{i(N(y_1^{(l)}+\cdots+y_k^{(l)})}-1} 
{e^{i((y_1^{(l)}+\cdots+y_k^{(l)})}-1} \prod_{j=1}^d 
\left(\prod_{s=k_1+\cdots+k_{j-1}+1}^{k_1+\cdots+k_j} Z_{G,j}(\,dy_s)\right),
\end{eqnarray*}
where we write $y=(y^{(1)},\dots,y^{(\nu)})$ for all 
$y\in[-\pi,\pi)^\nu$.  

The above sum of Wiener--It\^o integrals can be rewritten
with the change of variables $x_s=Ny_s$, $1\le s\le k$, in the 
form
\begin{eqnarray}
S_N&=& \sum_{\substack{(k_1,\dots,k_d),\; k_j\ge0,\;1\le j\le d,\\ 
k_1+\cdots +k_d=k}} \int c_{k_1,\dots,k_d} f^N(x_1+\cdots+x_k)
\nonumber \\
&&\qquad\qquad\qquad\qquad \prod_{j=1}^d 
\left(\prod_{s=k_1+\cdots+k_{j-1}+1}^{k_1+\cdots+k_j}
Z_{G^{(N)},j}(\,dx_s)\right),
\label{2.1}
\end{eqnarray}
where
\begin{equation}
f^N(x)=\prod_{l=1}^\nu  \frac{e^{i x^{(l)}}-1} 
{N(e^{ix^{(l)}/N}-1)} \label{2.2} 
\end{equation}
is a function on $[-N\pi,N\pi)^\nu$, and 
$Z_{G^{(N)},j}(A)=\frac{N^{\alpha/2}}{L(N)^{1/2}}Z_{G,j}(\frac AN)$ 
for all measurable sets $A\subset [-N\pi,N\pi)^\nu$ 
and $j=1,\dots,d$. In formula~(\ref{2.2}) the notation 
$x=(x^{(1)},\dots,x^{(\nu)})$ is applied for all
$x\in {\mathbb R}^\nu$. Let us observe that 
$(Z_{G^{(N)},1},\dots,Z_{G^{(N)},d})$ is a vector valued 
random spectral measure on the torus $[-N\pi,N\pi)^\nu$ 
which corresponds to the matrix valued spectral measure 
$G^{(N)}=(G^{(N)}_{j,j'})$, $1\le j,j'\le d$, on the 
torus $[-N\pi,N\pi)^\nu)$, defined by the formula 
$G^{(N)}_{j,j'}(A)=\frac{N^\alpha}{L(N)}G_{j,j'}(\frac AN)$, 
$1\le j,j'\le d$, on the sets 
$A\subset[-N\pi,N\pi)^\nu$, where 
$G=(G_{j,j'})$, $1\le j,j'\le d$, is the matrix valued 
spectral measure of the original vector valued stationary 
random field 
$X(p)=(X_1(p),\dots,X_d(p))$, $p\in{\mathbb Z}^\nu$.

In formulas (\ref{2.1}) and (\ref{2.2}) the normalized
random sum $S_N$ investigated in Theorem~3.2 is written
in the form of a sum of $k$-fold multiple Wiener--It\^o
integrals. Let us observe that the kernel functions 
$c_{k_1,\dots,k_d}f^N(x_1+\cdots+x_k)$ of these 
Wiener--It\^o integrals satisfy the relation
\begin{equation}
\lim_{N\to\infty} c_{k_1,\dots,k_d}f^N(x_1+\cdots+x_k)
=c_{k_1,\dots,k_d}f^0(x_1+\cdots+x_k) \label{2.3}
\end{equation}
for all indices $k_1,\dots, k_d$ such that $k_j\ge0$, 
$1\le j\le d$, and $k_1+\cdots+k_d=k$ with the function
\begin{equation}
f^0(x)=\prod_{l=1}^\nu \frac{e^{ix^{(l)}}-1} {ix^{(l)}} \label{2.4}
\end{equation}
defined on ${\mathbb R}^{\nu}$, and this convergence is 
uniform in all bounded subsets of ${\mathbb R}^{k\nu}$.

On the other hand, Proposition~3.1 states
that the elements of the matrix valued spectral measures
$G^{(N)}=(G^{(N)}_{j,j'})$ vaguely converge to the elements
of a matrix valued spectral measure $G^{(0)}=(G^{(0)}_{j,j'})$ on
${\mathbb R}^\nu$. In~(\ref{2.1}) we integrate with respect to
a vector valued random spectral measure corresponding to the
matrix valued spectral measure $(G^{(N)}_{j,j'})$, $1\le j,j'\le N$
of a generalized vector valued Gaussian stationary random field.
Hence it is natural to expect that the random variables $S_N$
converge in distribution to the random variable
\begin{eqnarray}
S_0&=&\sum_{\substack{(k_1,\dots,k_d),\; k_j\ge0,\;1\le j\le d,\\ 
k_1+\cdots +k_d=k}} 
\int c_{k_1,\dots,k_d}f^0(x_1+\dots+x_k) \nonumber \\
&&\qquad\qquad\qquad\qquad \prod _{j=1}^d 
\left(\prod_{s=k_1+\cdots+k_{j-1}+1}^{k_1+\cdots+k_j}
Z_{G^{(0)},j}(\,dx_s)\right), \label{2.5}
\end{eqnarray}
where $(Z_{G^{(0)},1},\dots,Z_{G^{(0)},d})$ is a vector valued random 
spectral measure on ${\mathbb R}^\nu$ corresponding to the matrix 
valued spectral measure $(G^{(0)}_{j,j'})$, $1\le j,j'\le d$. This 
is actually the statement of Theorem~3.2 with a slightly 
different notation.

\medskip
I shall formulate such a result in the following
Proposition~4A which helps to justify the above heuristic
argument. It states that this argument yields a correct
result if some additional conditions are also satisfied.
Theorem~3.2 will be proved with the help of this Proposition~4A
which is a reformulation of Proposition~3.1 in~\cite{11}.

\medskip
Before the presentation of Proposition~4A I recall from 
Section~5 of~\cite{10} the definition of that class of
functions which can be chosen for the kernel function of
a multiple Wiener--It\^o integral with respect to a
vector valued random spectral measure. This class of
functions appears in the formulation of~Proposition~4A.

\medskip
Let us consider the matrix valued spectral measure $G=(G_{j,j'})$,
$1\le j,j'\le d$, with non-atomic measures $G_{j,j}$,
$1\le j\le d$, of a vector valued Gaussian stationary
random field. (We can consider the spectral measure both
of an ordinary or of a generalized random field.) In
\cite{10} I have defined a real Hilbert space
${\cal K}_{k,j_1,\dots,j_k}={\cal K}_{k,j_1,\dots,j_k}
(G_{j_1,j_1},\dots,G_{j_k,j_k})$
depending on the diagonal elements $G_{1,1},\dots,G_{d,d}$
of the spectral measure~$G$ and on a sequence of integers
$(j_1,\dots,j_k)$ of length $k$ such that $1\le j_s\le d$ 
for all $1\le s\le k$. This Hilbert space has the property
that the $k$-fold Wiener--It\^o integral 
$$
\int f(x_1,\dots,x_k)Z_{G,{j_1}}(\,dx_1)\dots Z_{G,{j_k}}(\,dx_k)
$$
with respect to a vector valued random spectral measure
$$
Z_G=(Z_{G,1},\dots,Z_{G,d})
$$
corresponding to the matrix valued spectral measure $G$ is
defined for the kernel functions
$$
f(x_1,\dots,x_k)\in {\cal K}_{k,j_1,\dots,j_k}(G_{j_1,j_1},\dots,G_{j_k,j_k}).
$$
(In papers~\cite{10} and~\cite{11} I worked with Wiener--It\^o
integrals of order~$n$, while here I work with Wiener--It\^o
integrals of order~$k$. Hence I use here a slightly different
notation.)

We have $f\in{\cal K}_{k,j_1,\dots,j_k}(G_{j_1,j_1},\dots,G_{j_k,j_k})$
for a complex number valued function $f(x_1,\dots,x_k)$ with
arguments ${x_s\in\mathbb R}^{\nu}$, $1\le s\le k$,
if it satisfies the following conditions~(a) and~(b):

\medskip
\begin{description}
\item[{\rm(a)}] $f(-x_1,\dots,-x_k)=\overline{f(x_1,\dots,x_k)}$ 
for all $(x_1,\dots,x_k)\in {\mathbb R}^{k\nu}$,
\item[{\rm(b)}]
$\|f\|^2=\int|f(x_1,\dots,x_k)|^2
G_{j_1,j_1}(\,dx_1)\dots G_{j_k,j_k}(\,dx_n)<\infty$.
\end{description}

\medskip\noindent
The scalar product in
${\cal K}_{k,j_1,\dots,j_k}(G_{j_1,j_1},\dots,G_{j_k,j_k})$
is defined in the usual way. If 
$f,\,g\in{\cal K}_{k,j_1,\dots,j_k}(G_{j_1,j_1},\dots,G_{j_k,j_k})$,
then
$$
\langle f,g\rangle=\int f(x_1,\dots,x_k)\overline{g(x_1,\dots,x_k)}
G_{j_1,j_1}(\,dx_1)\dots G_{j_k,j_k}(\,dx_k).
$$

\medskip
In the formulation of Proposition~4A we take for all
$N=1,2,\dots$ a matrix valued non-atomic spectral
measure $G^{(N)}=(G^{(N)}_{j,j'})$,  $1\le j,j'\le d$, on
the torus $[-A_N\pi,A_N\pi)^\nu$ with a parameter~$A_N$
such that $A_N\to\infty$ as $N\to\infty$. We also take
 some functions 
$$
h^N_{j_1,\dots,j_k}(x_1,\dots,x_k)\in {\cal K}_{k,j_1,\dots,j_k}
={\cal K}_{k,j_1,\dots,j_k}(G^{(N)}_{j_1,j_1},\dots,G^{(N)}_{j_k,j_k})
$$
on the torus $[-A_N\pi,A_N\pi)^\nu$ for all sets of indices
$(j_1,\dots,j_k)$ such that $1\le j_s\le d$, 
$1\le s\le k$, and $N=1,2,\dots$. Besides, we fix for all 
$N=1,2,\dots$ a vector valued random spectral measure 
$Z_{G^{(N)}}=(Z_{G^{(N)},1},\dots,Z_{G^{(N)},d})$ on the torus 
$[-A_N\pi,A_N\pi)^\nu$ corresponding to the matrix valued 
spectral measure $G^{(N)}=(G^{(N)}_{j,j'})$, $1\le j,j'\le d$, and 
we define with the help of these quantities the  
sums of $k$-fold Wiener--It\^o integrals
\begin{equation}
S_N=\sum_{\substack{(j_1,\dots,j_k)\\ 
1\le j_s\le d, \textrm{ for all }1\le s\le k}}
\int h^N_{j_1,\dots,j_k}(x_1,\dots,x_k)Z_{G^{(N)},{j_1}}(\,dx_1)
\dots Z_{G^{(N)},{j_k}}(\,dx_k),    \label{2.7}
\end{equation}
$N=1,2,\dots$. We want to find some good conditions 
under which these random variables $S_N$ converge
in distribution to a random variable $S_0$, expressed
similarly as a sum of $k$-fold multiple
Wiener--It\^o integrals.

This will be done with the help of the following Proposition~4A which
agrees with Proposition~3.1 in paper~\cite{11}.

\medskip\noindent
{\bf Proposition~4A.} {\it Let us consider for all $N=1,2,\dots$
the sum of $k$-fold Wiener--It\^o integrals $S_N$ defined in 
formula~(\ref{2.7}) with the help of a vector valued random
spectral measure $Z_{G^{(N)}}=(Z_{G^{(N)},1},\dots,Z_{G^{(N)},d})$
corresponding to some non-atomic matrix valued spectral
measure $G^{(N)}=(G^{(N)}_{j,j'})$, $1\le j,j'\le d$, defined on 
a torus $[-A_N,A_N)^\nu$ such that $A_N\to\infty$ as 
$N\to\infty$ and functions 
$$
h^N_{j_1,\dots,j_k}(x_1,\dots,x_k)\in 
{\cal K}_{k,j_1,\dots,j_k}(G^{(N)}_{j_1,j_1},\dots,G^{(N)}_{j_k,j_k}).
$$
Let the coordinates $G^{(N)}_{j,j'}$, $1\le j,j'\le d$, of the 
matrix valued spectral measures $G^{(N)}=(G^{(N)}_{j,j'})$, $1\le j,j'\le d$, 
converge vaguely to the coordinates $G^{(0)}_{j,j'}$ of a non-atomic 
matrix valued spectral measure $G^{(0)}=(G^{(0)}_{j,j'})$, $1\le j,j'\le d$,
on ${\mathbb R}^\nu$ for all $1\le j,j'\le d$ as $N\to\infty$,
and let $Z_{G^{(0)}}=(Z_{G^{(0)},1},\dots,Z_{G^{(0)},d})$ be a vector valued
random spectral measure on ${\mathbb R}^\nu$ corresponding to the matrix
valued spectral measure $G^{(0)}=(G^{(0)}_{j,j'})$, $1\le j,j'\le d$. 
Let us also have some functions $h^0_{j_1,\dots,j_k}(x_1,\dots,x_k)$ on
${\mathbb R}^{k\nu}$ for all $1\le j_s\le d$, $1\le s\le k$,
such that these functions and matrix valued spectral measures 
satisfy the following conditions~(a) and~(b).

\medskip
\begin{description}
\item[\rm{(a)}] The functions $h^0_{j_1,\dots,j_k}(x_1,\dots,x_k)$ 
are continuous on ${\mathbb R}^{k\nu}$ for all $1\le j_s\le d$, 
$1\le s\le k$,  and for all $T>0$ and indices 
$1\le j_s\le d$, $1\le s\le k$, the functions 
$h^N_{j_1,\dots,j_k}(x_1,\dots,x_k)$ converge uniformly to 
the function $h^0_{j_1,\dots,j_k}(x_1,\dots,x_k)$ on the 
cube $[-T,T]^{k\nu}$ as $N\to\infty$.
\item[\rm{(b)}] For all $\varepsilon>0$ there is some 
$T_0=T_0(\varepsilon)>0$ such that 
\begin{equation}
\int_{{\mathbb R}^{k\nu}\setminus [-T,T]^{k\nu}}
|h^N_{j_1,\dots,j_k}(x_1,\dots,x_k)|^2
G^{(N)}_{j_1,j_1}(\,dx_1)\dots G^{(N)}_{j_k,j_k}(dx_k)<\varepsilon^2 
\label{2.8} 
\end{equation}
for all $1\le j_s\le d$, $1\le s\le k$, and $N=1,2\dots$ if $T>T_0$.
\end{description}

\medskip
Then
$$
h^0_{j_1,\dots,j_k}(x_1,\dots,x_k)\in{\cal K}_{k,j_1,\dots,j_k}
={\cal K}_{k,j_1,\dots,j_k}(G^{(0)}_{j_1,j_1},\dots G^{(0)}_{j_k,j_k}),
$$
inequality~(\ref{2.8}) holds also for $N=0$, the sum of
random integrals 
\begin{equation}
S_0=\sum_{\substack{(j_1,\dots,j_k)\\ 
1\le j_s\le d, \textrm{ for all }1\le s\le k}}
\int h^0_{j_1,\dots,j_k}(x_1,\dots,x_k)Z_{G^{(0)},{j_1}}(\,dx_1)
\dots Z_{G^{(0)},{j_k}}(\,dx_k) \label{2.9}
\end{equation}
exists, and the random variables $S_N$ defined in~(\ref{2.7}) 
converge to $S_0$ in distribution as $N\to\infty$.}

\medskip
(In the formulation of Proposition~4A I took the natural
identification of the torus $[-A_N,A_N)^\nu$ with the cube
$[-A_N,A_N)^\nu$ in the space ${\mathbb R}^\nu$. Thus, I considered
the functions $h_{j_1,\dots,j_k}(\cdot)$ as functions on
$\mathbb R^{k\nu}$ which disappear outside $[-A_N,A_n)^{k\nu}$,
and the complex measures $G_{j,j'}^{(N)}$ as complex measures on
${\mathbb R}^\nu$, concentrated on $[-A_N,A_N)^\nu$. In such a way
the vague convergence mentioned in the formulation of Proposition~4A
is meaningful.)  

\medskip
In the proof of Theorem~3.2 we want to show with the help of
Proposition~4A that the sequence of random variables
$S_N$, $N=1,2,\dots$, defined in~(\ref{2.1}) converge to the
random variable $S_0$ defined in~(\ref{2.5}) as $N\to\infty$.
To do this we rewrite these formulas with a different indexation
in such a way that the indices in the definition of the random
variables $S_N$ and $S_0$ fit to the indices in the definition
of the random variable of the random variables $S_N$ and $S_0$
appearing in the formulation of Proposition~4A. These random
variables were defined in formulas~(\ref{2.7}) and~(\ref{2.9}).

In formulas~(\ref{2.1}) and~(\ref{2.5}) summation is taken
for terms with indices \hfill\break
$(k_1,\dots,k_d)$ such that $k_s\ge0$,
$1\le s\le d$ and $k_1+\cdots+k_d=k$, while in the
corresponding expressions in formulas~(\ref{2.7})
and~(\ref{2.9}) in Proposition~4A it is taken for terms with
indices $(j_1,\dots,j_k)$ such that
$1\le j_s\le d$, $1\le s\le k$.

An important difference between the indexation in the two
cases is that in~(\ref{2.1}) and~(\ref{2.5}) only a special
subset of the indices in formulas~(\ref{2.7}) and~(\ref{2.9})
appear. Namely, if $s<s'$, and we compare the indices $j$
and $j'$ in the terms $Z_{G^{(N)},j}(\,dx_s)$ and
$Z_{G^{(N)},j'}(\,dx_{s'})$ belonging to these indices $s$ and
$s'$ in formula~(\ref{2.1}) or~(\ref{2.5}), then we find
that $j\le j'$. Hence such a reindexation of the indices
in~(\ref{2.1}) and~(\ref{2.5}) will be made, where the set
${\cal J}$ of the new indices is only a subset of the
indices $(j_1,\dots,j_k)$ appearing in formulas~(\ref{2.7})
and~(\ref{2.9}). Summation will be taken only for the
elements of~${\cal J}$ in these formulas.

More explicitly, the terms in the sums in~(\ref{2.1}) and~(\ref{2.5})
will be reindexed with such indices $(j_1,\dots,j_k)$ for which
the relation $1\le j_1\le j_2\le\cdots\le j_k\le d$ holds.
This is a subset of the set of indices $(j_1,\dots,j_k)$ appearing
in formulas~(\ref{2.7}) and~(\ref{2.9}). To carry out the desired
reindexation a one to one map will be defined between the sets
$$
{\cal J}=\{(j_1,\dots,j_k)\colon\;
1\le j_1\le j_2\le\cdots\le j_k\le d\}
$$
and 
$$
{\cal K}=\{(k_1,\dots,k_d)\colon\; k_s\ge0 \textrm{ for all }
1\le s\le d, \;\;k_1+\cdots+k_d=k\}.
$$

Put
\begin{eqnarray}
\textrm{for all }
(j_1,\dots,j_k)\in{\cal J}  
&&k_s(j_1,\dots j_k)=\textrm{the number of such elements } j_p 
\nonumber \\
&&\qquad \textrm{for which } j_p=s, \textrm{ for all } 1\le s\le d. 
\label{2.10}
\end{eqnarray}
This is a one to one map from ${\cal J}$ to ${\cal K}$ whose 
inverse is
\begin{eqnarray}
\textrm{for all } (k_1,\dots,k_d)\in{\cal K} 
&&j_s(k_1,\dots,k_d)=\min p\colon k_1+\cdots+k_p\ge s, \nonumber  \\ 
&&\qquad \textrm{for all } 1\le s\le k. \label{2.10a}
\end{eqnarray}
We shall apply these maps.

With the help of this correspondence between the sets ${\cal J}$
and ${\cal K}$ the random sums $S_N$ in~(\ref{2.1}) can be 
rewritten in a form where summation is taken for the sequences
$(j_1,\dots,j_k)\in{\cal J}$ instead of the sequences
$(k_1,\dots,k_d)\in{\cal K}$, and 
$k_s(j_1,\dots,j_k)$ is written instead of $k_s$, $1\le s\le d$.

The expression~$S_N$ defined in~(\ref{2.1}) can be rewritten as
\begin{equation}
S_N= \!\! \sum_{\substack{(j_1,\dots,j_k),\\ 
1\le j_1\le \cdots\le j_k\le d}}
\!\! \int c_{k_1(j_1,\dots,j_k),\dots,k_d(j_1,\dots,j_k)}f^N(x_1+\cdots+x_k) 
\prod_{s=1}^k Z_{G^{(N)},j_s}(\,dx_s)  \label{2.11}
\end{equation}
for all $N=1,2,\dots$ with the functions $f^N(x)$ defined 
in~(\ref{2.2}) and the indices $k_s(j_1,\dots,j_k)$, 
$1\le s\le d$ defined in~(\ref{2.10}).

To show why formula~(\ref{2.11}) holds let us rewrite
formula~(\ref{2.1}) (with the help of the one to one map we
defined between the sets ${\cal K}$ and ${\cal J}$) in the form
\begin{equation}
S_N=\!\!\! \sum_{\substack{(k_1,\dots,k_d),\; k_j\ge0,\;1\le j\le d,\\ 
k_1+\cdots +k_d=k}} \int c_{k_1,\dots,k_d} f^N(x_1+\cdots+x_k) 
\!\prod_{s=1}^k Z_{G^{(N)},j_s(k_1,\dots,k_d)}(\,dx_s). \label{2.12}
\end{equation}
To understand why formula~(\ref{2.12}) holds we have to show that
in (\ref{2.12}) the term $Z_{G^{(N)},u}(\,dx_s)$ with
$u=j_s(k_1,\dots,k_d)$ had to be chosen. It can be seen
from~(\ref{2.1}) that this number~$u$ must be chosen in such a way
that $k_1+\cdots+k_{u-1}+1\le s\le k_1+\cdots+k_u$. Then a comparison
of this condition with the definition of the mapping from~${\cal K}$
to ${\cal J}$ in~(\ref{2.10a}) shows that $u=j_s(k_1\dots,k_d)$.

Then if we rewrite the formula at the right-hand side of~(\ref{2.11})
by replacing the arguments $(k_1,\dots.k_d)\in{\cal K}$ by the arguments
$(j_1,\dots,j_k)\in{\cal J}$ with the help of the transformation we
defined from~${\cal J}$ to~${\cal K}$, and then we exploit that the
transformation we defined from~${\cal K}$ to~${\cal J}$ is its inverse
transformation, we get that formula~(\ref{2.12}) implies~(\ref{2.11}).

Relation~(\ref{2.11}) can be rewritten in the form
\begin{equation}
S_N=\sum_{\substack{(j_1,\dots,j_k),\\ 1\le j_1\le \cdots\le j_k\le d}}
 \int h^{N}_{j_1,\dots,j_k}(x_1,\dots,x_k) 
Z_{G^{(N)},{j_1}}(\,dx_1)\dots 
Z_{G^{(N)},{j_k}}(\,dx_k) \label{2.13}
\end{equation}
with
\begin{equation}
 h^{N}_{j_1,\dots,j_k}(x_1,\dots,x_k)= 
 c_{k_1(j_1,\dots,j_k),\dots,k_d(j_1,\dots,j_k)}f^N(x_1+\cdots+x_k),
 \label{2.14} 
\end{equation}
where the indices $k_s(j_1,\dots,j_k)$, $1\le s\le d$, are
defined in (\ref{2.10}).
Similarly, the random sum $S_0$ in~(\ref{2.5}) can 
be rewritten in the form
$$
S_0= \!\! \sum_{\substack{(j_1,\dots,j_k),\\ 
1\le j_1\le \cdots\le j_k\le d}}
\!\! \int c_{k_1(j_1,\dots,j_k),\dots,k_d(j_1,\dots,j_k)}f^0(x_1+\cdots+x_k) 
\prod_{s=1}^k Z_{G^{(0)},j_s}(\,dx_s)  
$$
with the function $f^0(x)$ defined
in~(\ref{2.4}) or in the following equivalent form.
\begin{equation}
S_0=\sum_{\substack{(j_1,\dots,j_k),\\ 1\le j_1\le \cdots\le j_k\le d}}
 \int h^{0}_{j_1,\dots,j_k}(x_1,\dots,x_k) 
Z_{G^{(0)},{j_1}}(\,dx_1)\dots 
Z_{G^{(0)},{j_k}}(\,dx_k) \label{2.15}
\end{equation}
with
\begin{equation}
 h^{0}_{j_1,\dots,j_k}(x_1,\dots,x_k)= 
 c_{k_1(j_1,\dots,j_k),\dots,k_d(j_1,\dots,j_k)}f^0(x_1+\dots+x_k).
 \label{2.16} 
\end{equation}

\section{Proof of the main theorems.}

Theorem~3.2 will be proved by means of the application
of Proposition~4A for the sequences  $S_N$ defined
in~(\ref{2.13}), (\ref{2.2}), (\ref{2.14})
and~(\ref{2.10}) for $N=1,2,\dots$, and in~(\ref{2.15}),
(\ref{2.4}), (\ref{2.16}) and~(\ref{2.10}) for $N=0$.
To do this we have to show that under the conditions
of Theorem~3.2 the conditions of Proposition~4A are also
satisfied with such a choice. Then the application of
Proposition~4A implies Theorem~3.2. (I would remark that the
random variable $S_0$ defined in formula~(\ref{1.8}) as the
limit in Theorem~3.2 agrees with the random variable~$S_0$
defined in~(\ref{2.5}), which is the same as the limit we get in
the application of Proposition~2A with the above written choice.
Only it is written there in a different form.)

To check the conditions of Proposition~4A let us first observe
that it follows from Proposition~3.1 that the (non-atomic)
elements $G^{(N)}_{j,j'}$ of the spectral measures $G^{(N)}$
vaguely converge to the (non-atomic) complex measures
$G^{(0)}_{j,j'}$ of a spectral measure~$G^{(0)}$ as $N\to\infty$
for all $1\le j,j'\le d$. It is also clear that the functions
$h^N_{j_1,\dots,j_k}(x_1,\dots,x_k)$ defined in~(\ref{2.14})
for all $1\le j_1\le\cdots \le j_k\le d$ and $N=1,2,\dots$ 
satisfy the condition $h^N_{j_1,\dots,j_k}(x_1,\dots,x_k)\in 
{\cal K}_{k,j_1,\dots,j_k}(G^{(N)}_{j_1,j_1},\dots,G^{(N)}_{j_k,j_k})$.

It follows from (\ref{2.3}), (\ref{2.14}) and~(\ref{2.16}) 
that condition~(a) of Proposition~4A holds with the functions
and measures chosen in the proof of  Theorem~3.2.  We
still have to prove relation~(\ref{2.8})  in condition~(b)
of Proposition~4A. This will be done with the help of the
following Proposition~5.1. (Actually in  Proposition~5.1
we prove a result slightly sharper than we need.)

\medskip\noindent
{\bf Proposition~5.1.} {\it Let us fix an integer $k\ge1$,
and let $G=(G_{j,j'})$, $1\le j,j'\le d$, be the matrix valued 
spectral measure of a vector valued stationary random field 
$X(p)=(X_1(p),\dots,X_p(d))$, $p\in{\mathbb Z}^\nu$, 
defined on the torus $[-\pi,\pi)^\nu$ with such correlation 
function $r_{j,j'}(p)=EX_j(0)X_{j'}(p)$, $1\le j,j'\le d$, 
$p\in{\mathbb Z}^\nu$, that satisfies relation~(\ref{1.3}) with
some $0<\alpha<\frac\nu k$. For all $N=1,2,\dots$ let us
consider the measures $G^{(N)}_{j,j}$, $1\le j\le d$, defined
in formula~(\ref{1.6}) together with the measures 
$\mu^{(N)}_{j_1,\dots,j_k}$ defined for all sets of indices
$j_1,\dots,j_k$ such that $1\le j_s\le d$, $1\le s\le k$, on 
${\mathbb R}^{k\nu}$ by the formula
\begin{equation}
\mu^{(N)}_{j_1,\dots,j_k}(A)=\int_A|h_N(x_1,\dots,x_k)|^2 
G^{(N)}_{j_1,j_1}(\,dx_1)\dots G^{(N)}_{j_k,j_k}(\,dx_k), 
\quad A\in {\cal B}^{k\nu}, \label{2.17}
\end{equation}
with
\begin{equation}
h_N(x_1,\dots,x_k)=f^N(x_1+\cdots+x_k)
=\prod_{l=1}^\nu
\frac{e^{i((x_1^{(l)}+\cdots+x_k^{(l)})}-1} 
{N(e^{i((x_1^{(l)}+\cdots+x_k^{(l)})/N}-1)}, \label{2.18} 
\end{equation}
where we use the notation $x=(x^{(1)},\dots,x^{(\nu)})$ 
for a vector $x\in {\mathbb R}^\nu$. These measures 
$\mu^{(N)}_{j_1,\dots,j_k}$ converge weakly to a finite 
measure $\mu^{(0)}_{j_1,\dots,jk}$ on ${\mathbb R}^{k\nu}$.} 

\medskip\noindent
{\it Proof of Theorem~3.2 with the help of Proposition~5.1.}
As we have seen to prove Theorem~3.2 it is enough to check that
the measures $G^{(N)}=(G^{(N)}_{j,j'})$, $1\le j,j'\le d$, and
functions $h^N_{j_1,\dots,j_k}$ defined before satisfy the conditions
of Proposition~4A, since this enables us to apply this result.
Moreover, we have proved the validity of all of these conditions
except formula~(\ref{2.8}) in condition~(b) of Proposition~4A.
But the validity of this condition follows from Proposition~5.1,
since this result implies that the measures $\mu^{(N)}_{j_1,\dots,j_k}$,
$N=1,2,\dots$, defined in~(\ref{2.17}) and~(\ref{2.18}) are
uniformly tight. This fact together with the definition of the
measures $\mu^{(N)}_{j_1,\dots,j_k}$ and the identity
$h^N_{j_1,\dots,j_k}(x_1,\dots,x_k)=
c_{k_1(j_1,\dots,j_k),\dots,k_d(j_1,\dots,j_k)} h_N(x_1,\dots,x_k)$
imply that relation~(\ref{2.8}) holds. Theorem~3.2 is proved. 

\medskip
It remained to prove Proposition~5.1.

\medskip\noindent
{\it Proof of Proposition~5.1.} Most calculations needed
in the proof of Proposition~5.1 were actually carried out
in the proof of Theorem~8.2 of~\cite{9}.  Only some slight
modifications are needed in the proof. In some steps I
shall refer to the corresponding part in the proof of
Theorem~8.2 in~\cite{9} and omit the details of the
calculation. 

I compute for all $N=1,2,\dots$ the Fourier transform
$$
\varphi^{(N)}_{j_1,\dots,j_d}(t_1,\dots,t_k)=\int
e^{i((t_1,x_1)+\cdots+(t_k,x_k))}
\mu^{(N)}_{j_1,\dots,j_k}(\,dx_1,\dots,dx_k)
$$
of the measures $\mu^{(N)}_{j_1,\dots,j_k}$ defined in~(\ref{2.17}) 
and give a good asymptotic formula for it. More precisely
I do this only for such coordinates $(t_1,\dots,t_k)$ of 
the function $\varphi^{(N)}_{j_1,\dots,j_d}(t_1,\dots,t_k)$ 
which have the form $t_l=\frac {p_l}N$ with some 
$p_l\in {\mathbb Z}^\nu$, $l=1,\dots,k$. But as it is explained
at the end of this proof, even such a result is sufficient for
us. In the calculation of the formula expressing 
$\varphi^{(N)}_{j_1,\dots,j_d}(t_1,\dots,t_k)$ I exploit that
the function $h_N(x_1,\dots,x_k)$ defined in~(\ref{2.18}) can be 
written in the form
$$
h_N(x_1,\dots,x_k)=\frac1{N^\nu}\sum_{u\in B_N}
\exp\left\{i\frac1N(u,x_1+\cdots+x_k)\right\}.
$$
Hence, and because of the definition of the spectral measures
$G^{(N)}_{j,j}(\cdot)$ in~(\ref{1.6})
\begin{eqnarray*}
&&\!\!\!\!\!\!\!\!\!
\varphi^{(N)}_{j_1,\dots,j_d}(t_1,\dots,t_k)
=\frac1{N^{2\nu}}\int\exp\left\{
i\frac1N((p_1,x_1)+\cdots+(p_k,x_k))\right\}    \\
&&\sum_{u\in B_N}\sum_{v\in B_N}
\exp\left\{i\left(\frac{u-v}N,x_1+\cdots+x_k\right)\right\} 
G^{(N)}_{j_1,j_1}(\,dx_1)\dots G^{(N)}_{j_k,j_k}(\,dx_k)   \\
&&\qquad=\frac1{N^{2\nu}} \sum_{u\in B_N}\sum_{v\in B_N}\left(\prod_{s=1}^k
\int \exp\left\{i\left(\frac{u-v+p_s}N ,x_s\right)\right\}
 G^{(N)}_{j_s,j_s}(\,dx_s)\right) \\
&&\qquad=\frac1{N^{2\nu-k\alpha}L(N)^k}\sum_{u\in B_N}\sum_{v\in B_N}
r_{j_1,j_1}(u-v+p_1)\cdots r_{j_k,j_k}(u-v+p_k)
\end{eqnarray*}
if $t_l=\frac {p_l}N$ with some $p_l\in{\mathbb Z}^\nu$, 
$1\le l\le k$. This identity can be rewritten by taking 
the summation at the right-hand side of the last formula 
first for such pairs $(u,v)$ for which $u-v=y$ with 
a fixed point $y\in{\mathbb Z}^\nu$ and then for 
the lattice points $y\in{\mathbb Z}^\nu$. By working
with $x=\frac yN$ instead of $y$ we get that  
$$
\varphi^{(N)}_{j_1,\dots,j_d}(t_1,\dots,t_k)
=\int_{[-1,1]^\nu} f^{(N)}_{j_1,\dots,j_k}(t_1,\dots,t_k,x)\lambda_N(\,dx)
$$
with
\begin{eqnarray*}
&& \!\!\!\!\!\!\!\!\!\!\!\!\!
f^{(N)}_{j_1,\dots,j_k}(t_1,\dots,t_k,x) \\
&&\!\!\!\!\!  =\left(1-\frac{|x^{(1)}N|}N\right)\cdots
\left(1-\frac{|x^{(\nu)}N|}N\right)
\frac{r_{j_1,j_1}(N(x+t_1))}{N^{-\alpha}L(N)}\cdots
\frac{r_{j_k,j_k}(N(x+t_k))}{N^{-\alpha}L(N)},
\end{eqnarray*}
where $\lambda_N$ is the measure concentrated in the points of 
the form $x=\frac pN$ with such points
$p=(p_1,\dots,p_\nu)\in {\mathbb Z}^\nu$ for which $-N<p_l<N$
for all $1\le l\le\nu$, and $\lambda_N(x)=N^{-\nu}$ for each 
point~$x$ with this property. (Here such a calculation is
applied which is similar to that in the proof of Theorem~8.2
of~\cite{9} when formula~(8.20) of that work was rewritten
in another form.)

Let us extend the definition of
$\varphi^{(N)}_{j_1,\dots,j_d}(t_1,\dots,t_k)$ to all
$(t_1,\dots,t_k)\in{\mathbb R}^{k\nu}$ by defining it as
$$
\varphi^{(N)}_{j_1,\dots,j_d}(t_1,\dots,t_k)
=\varphi^{(N)}_{j_1,\dots,j_d}\left(\frac{p_1}N,\dots,\frac{p_k}N\right),
\quad t_l\in{\mathbb R}^\nu \textrm{ for all } 1\le l\le k,
$$
where $p_l=p_l(t_l)$ is defined as the integer part $[t_lN]$ of $t_lN$,
$1\le l\le k$, i.e. $p_l\in{\mathbb Z}^\nu$, and 
$p_l^{(s)}\le t_l^{(s)}N< p_l^{(s)}+1$ if $t_l^{(s)}>0$, and
$p_l^{(s)}-1<t_l^{(s)}N\le p_l^{(s)}$ if $t_l^{(s)}\le0$, $1\le s\le\nu$.

Let us also extend the definition of the function
$f^{(N)}_{j_1,\dots,j_k}(t_1,\dots,t_k,x)$ to 
$(t_1,\dots,t_k,x)\in{\mathbb R}^{k\nu}\times[-1,1]^\nu$ by means 
of the formula
\begin{eqnarray*}
&& \!\!\!\!\!\!\!\!\!\!\!\!\!
f^{(N)}_{j_1,\dots,j_k}(t_1,\dots,t_k,x) \\ 
&&\!\!\!\!\! =\left(1-\frac{|q^{(1)}|}N\right)\cdots
\left(1-\frac{|q^{(\nu)}|}N\right)
\frac{r_{j_1,j_1}(q+p_1)}{N^{-\alpha}L(N)}\cdots
\frac{r_{j_k,j_k}(q+p_k)}{N^{-\alpha}L(N)}
\end{eqnarray*}
for $t_l\in{\mathbb R}^\nu$, $1\le l\le k$, and $x\in[-1,1]^\nu$,
where $p_l=p_l(t)$ is defined as before, and $q=q(x)$ is defined
as $q=(q^{(1)},\dots,q^{(\nu)})\in{\mathbb Z}^\nu$
with $q_l\in{\mathbb Z}^\nu$, and 
$q^{(s)}\le x_l^{(s)}< q^{(s)}+1$.

We have
\begin{equation}
\varphi^{(N)}_{j_1,\dots,j_d}(t_1,\dots,t_k)
=\int_{[-1,1]^\nu} f^{(N)}_{j_1,\dots,j_k}(t_1,\dots,t_k,x)\,dx \label{2.18a}
\end{equation}
for the functions $\varphi^{(N)}(\cdot)$ and $f^{(N)}(\cdot)$
with this extended domain of definition, where $dx$ denotes
integration with respect to the Lebesgue measure.

It follows from relation~(\ref{1.3}) and the fact that $\frac qN$
is very close to $x$, and $\frac{p_l}N$ is very close to $t_l$,
for all $1\le l\le k$ if $N$ is large that for all parameters
$t_1,\dots,t_k$ and $\varepsilon>0$
$$
f^{(N)}_{j_1,\dots,j_k}(t_1,\dots,t_k,x)
\to f^{(0)}_{j_1,\dots,j_k}(t_1,\dots,t_k,x) 
$$
holds uniformly with the limit function
\begin{eqnarray*}
&&f^{(0)}_{j_1,\dots,j_k}(t_1,\dots,t_k,x) \\
&&\qquad =(1-|x^{(1)}|)\dots(1-|x^{(\nu)}|)
\frac{a_{j_1,j_1}\left(\frac{x+t_1}{|x+t_1|}\right)}{|x+t_1|^\alpha}\dots
\frac{a_{j_k,j_k}\left(\frac{x+t_k}{|x+t_k|}\right)}{|x+t_k|^\alpha}
\end{eqnarray*}
on the set $x\in[-1,1]^\nu \setminus
\bigcup\limits_{l=1}^k\{x\colon\;|x+t_l|<\varepsilon\}$. 

Some additional calculation shows that for small $\varepsilon>0$
integration on the domain
\begin{eqnarray*}
&&[-1,1]^\nu \setminus \left([-1,1]^\nu \setminus
\bigcup\limits_{l=1}^k\{x\colon\;|x+t_l|<\varepsilon\}\right) \\
&&\qquad\qquad\qquad  =
[-1,1]^\nu\cap\left(\bigcup_{l=1}^k \{x\colon\; |x+t_l|<\varepsilon\}\right)
\end{eqnarray*}
gives a negligible contribution to the integral in
formula~(\ref{2.18a}) (with parameters $j_1,\dots,j_k$ and
$t_1,\dots,t_l$), or to the integral that we get if the kernel
function $f^{(N)}_{j_1,\dots,j_k}$ is replaced by $f^{(0)}_{j_1,\dots,j_k}$
in the integral in~(\ref{2.18a}). Hence the relation
\begin{equation}
\varphi^{(N)}_{j_1,\dots,j_k}(t_1,\dots,t_k)\to
\varphi^{(0)}_{j_1,\dots,j_k}(t_1,\dots,t_k)
=\int_{[-1,1]^\nu}f^{(0)}_{j_1,\dots,j_k}(t_1,\dots,t_k,x)\,dx
\label{2.19}
\end{equation}
holds for all $(t_1,\dots,t_k)\in {\mathbb R}^{k\nu}$ as $N\to\infty$,
and $\varphi^{(0)}_{j_1,\dots,j_k}(t_1,\dots,t_k)$ is a continuous function.
This calculation was carried out in that part of the proof of
Theorem~8.2 in~\cite{9} which followed the discussion of Lemma~8.4.
Hence here I omit it.

By a classical result of probability theory if the Fourier 
transforms of a sequence of finite measures on ${\mathbb R}^{k\nu}$
converge to a function continuous at the origin, then the
limit function is also the Fourier transform of a finite
measure on ${\mathbb R}^{k\nu}$, and the sequence of probability 
measures whose Fourier transforms were taken converge to 
this measure. In the proof of Proposition~5.1 this result
cannot be applied, because we have a control on the Fourier 
transform of $\mu^{(N)}_{j_1,\dots,j_k}$ only in points of 
the form $(t_1,\dots,t_k)$ with $t_l=\frac{p_l}N$ and
$p_l\in{\mathbb Z}^\nu$, $1\le l\le k$. But the measures
$\mu^{(N)}_{j_1,\dots,j_k}$ have the additional property that
they are concentrated in the cube $[-N\pi,N\pi)^{k\nu}$.
Lemma~8.4 of~\cite{9} can be applied, and it shows that
relation~(\ref{2.19}) and the continuity of the limit function 
$\varphi^{(0)}_{j_1,\dots,j_k}(t_1,\dots,t_k)$ together with
the above mentioned concentration property of the measures
$\mu^{(N)}_{j_1,\dots,j_k}$ imply the weak convergence of the
measures~$\mu^{(N)}_{j_1,\dots,j_k}$ to a finite measure
$\mu^{(0)}_{j_1,\dots,j_k}$. This result also implies that
this finite measure $\mu^{(0)}_{j_1,\dots,j_k}$ has the Fourier
transform $\varphi^{(0)}_{j_1,\dots,j_k}(t_1,\dots,t_k)$.
Proposition~5.1 is proved.

\medskip
To prove Theorem~3.3 with the help of Theorem~3.2 it is
enough to show that if a function $H^{(1)}(\cdot)$ satisfies
(\ref{1.9}) and (\ref{1.10}), and the Gaussian stationary
random field $X(p)=(X_1(p),\dots,X_d(p))$ satisfies~(\ref{1.3})
and~(\ref{1.4}), then 
\begin{equation}
\frac1{N^{\nu-k\alpha/2}L(N)^{k/2}}\sum_{p\in B_N}H^{(1)}(X_1(p),\dots,X_d(p))
\Rightarrow0 \quad \textrm{as }N\to\infty, \label{2.20}
\end{equation}
where $\Rightarrow$ denotes convergence in probability. I 
shall prove that even the second moments of the normalized 
sums in~(\ref{2.20}) tend to zero as $N\to\infty$. The following
Lemma~5A which agrees with Lemma~1 of~\cite{1} (only with a
slightly different notation) helps in the proof of this statement.

\medskip\noindent
{\bf Lemma~5A.} {\it Let $X=(X_1,\dots,X_d)$ and $Y=(Y_1,\dots Y_d)$
be two Gaussian random vectors with expectation zero 
such that $EX_jX_{j'}=EY_jY_{j'}=\delta_{j,j'}$, $1\le j,j'\le d$, 
and let $r_{j,j'}=EX_jY_{j'}$, $1\le j,j'\le d$. Take a number 
$k\ge1$ and a function $H^{(1)}(x_1,\dots,x_d)$ that satisfies 
relations~(\ref{1.9}) and~(\ref{1.10}). Assume that
$$
\psi=\max\left( \left(\sup_{1\le j\le d}\sum_{j'=1}^d |r_{j,j'}|\right),
\left(\sup_{1\le j'\le d}\sum_{j=1}^d |r_{j,j'}|\right)\right)\le 1.
$$
Then
$$
|EH^{(1)}(X_1,\dots,X_d)H^{(1)}(Y_1,\dots,Y_d)|\le \psi^{k+1}
E\left[{H^{(1)}}(X_1,\dots,X_d)\right]^2.
$$
}

\medskip\noindent
{\it Proof of Theorem~3.3.} It follows from relations (\ref{1.3}),
(\ref{1.4}) and Lemma~5A together with the inequality
$E\left[H^{(1)}(X_1(0),\dots,X_d(0))\right]^2<\infty$ which holds
because of~(\ref{1.10}) that for two elements
$X(p)=(X_1(p),\dots,X_d(p))$ and $X(q)=(X_1(q),\dots,X_d(q))$,
$p,q\in{\mathbb Z}^\nu$, of our vector valued Gaussian stationary
random field there exists some threshold index $n_0\ge1$ and
constant $0<C<\infty$ such that 
\begin{eqnarray*}
&&|EH^{(1)}(X_1(p),\dots X_d(p))H^{(1)}(X_1(q),\dots,X_d(q))| \\
&&\qquad\qquad \le C|p-q|^{-(k+1)\alpha}L(|p-q|)^{k+1}
\end{eqnarray*}
if $|p-q|\ge n_0$. On the other hand,
\begin{eqnarray*}
&&|EH^{(1)}(X_1(p),\dots X_d(p))H^{(1)}(X_1(q),\dots,X_d(q))|\\
&&\qquad \qquad \le E{H^{(1)}}^2(X_1(0),\dots X_d(0))\le C_1
\end{eqnarray*}
for all $p,q\in{\mathbb Z}^\nu$ with some $C_1<\infty$ by the 
Schwarz inequality and relation~(\ref{1.10}). Hence we get by summing
up the above two inequalities for all $q\in B_N$ with a fixed $p\in B_N$,
and applying the first inequality if $|p-q|>n_0$ and the second one if
$|p-q|\le n_0$ that
\begin{eqnarray*}
&&\left|EH^{(1)}(X_1(p),\dots,X_d(p))
\left(\sum_{q\in B_N}H^{(1)}(X_1(q),\dots,X_d(q)\right)\right| \\
&&\qquad\qquad\qquad \le C_2(1+N^{\nu+\varepsilon-(k+1)\alpha})
\end{eqnarray*}
for all $p\in B_N$ and $\varepsilon>0$ with an appropriate 
$C_2=C_2(\varepsilon)>0$. Since $\nu-k\alpha>0$ we get by
summing up the last inequality for all $p\in B_N$ that
$$
\frac1{N^{2\nu-k\alpha}L(N)^k}
E\left[\sum_{p\in B_N}H^{(1)}(X_1(p),\dots,X_d(p))\right]^2\to0 
\textrm\quad \textrm{as }N\to\infty.
$$
Indeed, it can be seen that for all $\varepsilon>0$ the
expression in the last formula can be bounded from
above by $C(\varepsilon)N^{-\delta+\varepsilon}$ with
$\delta=\min(\nu-k\alpha,\alpha)>0$ and a constant
$C(\varepsilon)>0$ depending only on $\varepsilon$. This
implies formula~(\ref{2.20}). Formula~(\ref{2.20}) together
Theorem~3.2 yield Theorem~3.3. Theorem~3.3 is proved.

\medskip
{\it Proof of Theorem~3.4.}\/  The proof of Theorem~3.4 is
very similar to that of Theorems~3.2 and~3.3. Hence I only
briefly explain it.

It is enough to show that for any positive integer $K$, 
parameters $t_1,\dots,t_K$, $t_p\in[0,1]^\nu$, 
$1\le p\le K$ and real constants $C_1,\dots,C_K$ the linear 
combinations $\sum_{p=1}^K C_p S_N(t_p)$ converge to
$\sum_{p=1}^K C_p S_0(t_p)$ in distribution as $N\to\infty$,
since this implies that the random vectors 
$(S_N(t_1),\dots,S_N(t_K))$ converge in distribution to 
the random vector $(S_0(t_1),\dots,S_0(t_K))$ as 
$N\to\infty$. Moreover, similarly to the proof of 
Theorem~3.3 the proof of Theorem~3.4 can be reduced
to the case $H(x_1,\dots,x_d)=H^{(0)}(x_1,\dots,x_d)$ 
with a function $H^{(0)}(x_1,\dots,x_d)$ which satisfies
relation~(\ref{1.5}).

In the first step of the proof the linear combinations
$\sum_{p=1}^K C_p S_N(t_p)$, $N=0,1,2,\dots$, are written 
in the form of a sum of $k$-fold Wiener-It\^o integrals 
with respect to the coordinates of an appropriate vector 
valued random spectral measure. This can be done, first
by writing the random variables $S_N(t)$ for all
$t\in[0,1]^\nu$ in the desired form. The random
variables $S_0(t)$ are written in such a form
in~(\ref{1.13}). In the case $N=1,2,\dots$ the right 
expression of $S_N(t)$ in the form of a sum of
Wiener--It\^o integrals can be found similarly to the
method applied in the proof of Theorem~3.2. We can
write, similarly to the proof of formulas~(\ref{2.1})
and~(\ref{2.2})  
\begin{eqnarray*}
S_N(t)&=& 
\sum_{\substack{(k_1,\dots,k_d),\; k_j\ge0,\;1\le j\le d,\\ 
k_1+\cdots +k_d=k}}\int c_{k_1,\dots,k_d} 
f^N(t,x_1+\cdots+x_k) \\
&&\qquad\qquad\qquad \prod_{j=1}^d 
\left(\prod_{l=k_1+\cdots+k_{j-1}+1}^{k_1+\cdots+k_j} 
Z_{G^{(N)},j}(\,dx_l)\right)
\end{eqnarray*}
with
$$
f^N(t,x)=\prod_{l=1}^\nu 
\frac{\exp\left\{i\frac{]t^{(l)}N[}N(x^{(l)})\right\}-1} 
{N\left(\exp\left\{i\frac1N(x^{(l)})\right\}-1\right)}, 
$$
where $t=(t^{(1)},\dots,t^{(\nu)})$, the number $]t^{(l)}N[$
in the definition of the function \hfill\break 
$f^N(t,x_1,\dots,x_k)$) is the smallest integer which is
not smaller than $t^{(l)}N$, and $Z_{G^{(N)},j}$ agrees with
the spectral measure that appeared in formula~(\ref{2.1}).

It is not difficult to see that, similarly to relations 
(\ref{2.3}) and (\ref{2.4})
$$
\lim_{N\to\infty} f^N(t,x_1+\cdots+x_k)=f^0(t,x_1+\cdots+x_k) 
$$
with the function
$$
f^0(t,x)=\prod_{l=1}^\nu \frac{e^{it^{(l)}(x^{(l)})}-1} {i(x^{(l)})} 
$$
for all $(x_1,\dots,x_k)\in{\mathbb R}^{k\nu}$, and for a fixed
parameter~$t$ this convergence is uniform in all bounded
subsets of ${\mathbb R}^{k\nu}$.

With the help of the above considerations the proof of 
Theorem~3.4 can be reduced, similarly to the proof of 
Theorem~3.2 to the following statement. 

Fix some number~$K$, real constants $C_1,\dots,C_K$ and 
points $t_1,\dots t_K$ with $t_p\in [0,1]^\nu$, $1\le p\le K$ 
together with some constants $c_{k_1,\dots,k_d}$ with parameters 
$k_j\ge0$, $1\le j\le d$, and $k_1+\cdots+k_d=k$ which agree
with the coefficients in the sum~(\ref{1.5}). Let us define
with their help the random sums
\begin{eqnarray}
S_N&=& 
\sum_{\substack{(k_1,\dots,k_d),\; k_j\ge0,\;1\le j\le d,\\ 
k_1+\cdots +k_d=k}} \int \left(\sum_{p=1}^K C_p c_{k_1,\dots,k_d} 
f^N(t_p,x_1+\cdots+x_k)\right) \nonumber \\ 
&&\qquad\qquad\qquad \prod_{j=1}^d 
\left(\prod_{l=k_1+\cdots+k_{j-1}+1}^{k_1+\cdots+k_j} 
Z_{G^{(N)},j}(\,dx_l)\right) \label{2.21} 
\end{eqnarray}
with the above defined functions $f^{(N)}(t,x_1,\dots,x_k)$
for all $N=1,2,\dots$, and
\begin{eqnarray}
S_0&=& 
\sum_{\substack{(k_1,\dots,k_d),\; k_j\ge0,\;1\le j\le d,\\ 
k_1+\cdots +k_d=k}} 
\int \left(\sum_{p=1}^KC_p c_{k_1,\dots,k_d} 
f^0(t_p,x_1+\cdots+x_k)\right) \nonumber \\ 
&&\qquad\qquad\qquad \prod_{j=1}^d 
\left(\prod_{l=k_1+\cdots+k_{j-1}+1}^{k_1+\cdots+k_j} 
Z_{G^{(0)},j}(\,dx_l)\right) \label{2.22} 
\end{eqnarray}
with the previously defined function 
$f^0(t,x_1,\dots,x_k)$. The sequence of 
random variables $S_N$ defined in~(\ref{2.21}) converge in 
distribution to $S_0$ defined in~(\ref{2.22}) as $N\to\infty$.

\medskip
This statement can be proved, similarly to Theorem~3.2 
with the help of Proposition~4A. First the random variables
$S_N$, $N=1,2,\dots$, and $S_0$ must be rewritten in a form
in which Proposition~4A can be applied. They can be
rewritten in the form of a sum of multiple Wiener--It\^o
integrals indexed by sequences of integers $j_1,\dots,j_k$
such that $1\le j_1\le\cdots\le j_k\le d$. This can be
done similarly to the rewriting of formulas~(\ref{2.1}) 
and~(\ref{2.5}) in formulas~(\ref{2.13}), (\ref{2.14}) 
and~(\ref{2.15}), (\ref{2.16}) with the help of the
expressions $k_s(j_1\dots,j_k)$ defined in~(\ref{2.10}).
The random variable $S_N$ in (\ref{2.21}) can be
rewritten as
\begin{eqnarray}
S_N  \!\! &=&  \!\! \sum_{\substack{(j_1,\dots,j_k),\\ 
1\le j_1\le\cdots\le j_k\le d}} \!\! \int \left(\sum_{p=1}^K
C_p c_{k_1(j_1,\dots,j_k),\dots,k_d(j_1,\dots,j_k)}
f^N(t_p,x_1+\cdots+x_k)\right) \nonumber \\
&&\qquad\qquad\qquad\qquad \quad
Z_{G^{(N)},{j_1}}(\,dx_1)\dots Z_{G^{(N)},{j_k}}(\,dx_k) \label{2.23}
\end{eqnarray}
for all $N=1,2,\dots$, and the random variable in~(\ref{2.22}) as
\begin{eqnarray}
S_0&=& \sum_{\substack{(j_1,\dots,j_k),\\
1\le j_1\le \cdots\le j_k\le d}} \int\left(\sum_{p=1}^K
C_p c_{k_1(j_1,\dots,j_k),\dots,k_d(j_,\dots,j_k)} 
f^0(t_p,x_1+\cdots+x_k)\right) \nonumber \\
&&\qquad\qquad\qquad\qquad \quad
Z_{G^{(0)},{j_1}}(\,dx_1)\dots Z_{G^{(0)},{j_k}}(\,dx_k), 
\label{2.24}
\end{eqnarray}
where the indices $k_s(j_1,\dots,j_k)$, $1\le s\le d$, are 
defined in~(\ref{2.10}).

The random integrals in formulas~(\ref{2.23}) and~(\ref{2.24})
have kernel functions of the form
\begin{eqnarray}
&&h^{N}_{j_1,\dots,j_k}(x_1,\dots,x_k)
=h^{N}_{j_1\dots,j_k,t_1,\dots,t_K}(x_1,\dots,x_k) \label{2.25}\\
&&\qquad= \sum_{p=1}^K C_p c_{k_1(j_1,\dots,j_k),\dots,k_d(j_1,\dots,j_k)}
f^N(t_p,x_1+\cdots+x_k) \nonumber
\end{eqnarray}
for all $N=0,1,2,\dots$. Let us define for all $N=0,1,2,\dots$
the measures $\mu_{N,j_1,\dots,j_k}$ as
\begin{eqnarray}
\mu^{(N)}_{j_1,\dots,j_k}(A)&=&
\mu^{(N)}_{j_1,\dots,j_k,t_1,\dots,t_K}(A)
\label{2.26} \\
&=&\int_A|h^N_{j_1,\dots,j_k,t_1,\dots,t_K}(x_1,\dots,x_k)|^2  
G^{(N)}_{j_1,j_1}(\,dx_1)\dots G^{(N)}_{j_k,j_k}(\,dx_k) 
\nonumber
\end{eqnarray}
where integral is taken for all measurable sets $A\in{\cal B}^{k\nu}$.

We want to show with the help of Proposition~4A that the
distributions of the random variables $S_N$, $N=1,2,\dots$,
defined in~(\ref{2.23}) converge weakly to the distribution of
the random variable $S_0$ defined in~(\ref{2.24}).
This implies Theorem~3.4.

To prove this convergence we have to show that the functions
$h^N_{j_1,\dots,j_k}$, $N=0,1,2,\dots$, defined in (\ref{2.25}) 
and the measures $G^{(N)}_{j,j}$, $1\le j\le d$, 
$N=0,1,2,\dots$, satisfy the conditions of Proposition~4A. 
The main point is to prove relation~(\ref{2.8}) in 
condition~(b) of Proposition~4A. To prove this we show 
that the measures $\mu^{(N)}_{j_1,\dots,j_k}$, $N=1,2,\dots$, 
defined in (\ref{2.26}) are tight, i.e. for all $\varepsilon>0$ 
there exists a $T=T(\varepsilon,j_1,\dots,j_k,t_1,\dots,t_K)$ 
such that 
$$
\mu^{N}_{j_1,\dots,j_k,t_1,\dots,t_K}
(\mathbb R^{k\nu}\setminus[-T,T]^{k\nu})<\varepsilon \quad
\textrm{for all }N=1,2,\dots.
$$

Because of the Schwarz inequality and the definition 
of the functions \hfill\break 
$h^{N}_{j_1\dots,j_k,t_1,\dots,t_K}(x_1,\dots,x_k)$
the proof of this tightness property can be reduced to the
justification of the following inequality.

Let us define for all $t=(t_1,\dots,t_\nu)\in[0,1]^\nu$,
and $N=1,2,\dots$ the measure $\mu_{N,t}$ on 
${\mathbb R}^{k\nu}$ by the formula
\begin{eqnarray*}
\mu_{N,t}(A)&=& \int_A|f^N(t,x_1+\cdots+x_k)|^2
G^{(N)}_{j_1,j_1}(\,dx_1)\dots G^{(N)}_{j_k,j_k}(\,dx_k)\\ 
&=&\int_A \left|\prod_{l=1}^\nu 
\frac{\exp\left\{i\frac{]t^{(l)}N[}N(x_1^{(l)}+\cdots+x_k^{(l)})\right\}-1} 
{N\left(\exp\left\{i\frac1N(x_1^{(l)}+\cdots+x_k^{(l)})\right\}-1
\right)}\right|^2 \\
&&\qquad\qquad \qquad\qquad \qquad
G^{(N)}_{j_1,j_1}(\,dx_1)\dots G^{(N)}_{j_k,j_k}(\,dx_k) 
\end{eqnarray*}
for all $A\in{\cal B}^{k\nu}$. The inequality
$$
\mu_{N,t}({\mathbb R}^{k\nu}\setminus [-T,T]^{k\nu})<\varepsilon
$$
holds for all $N=1,2,\dots$, if $T\ge T_0(\varepsilon,t)$ with an
appropriate threshold index $T_0(\varepsilon,t)>0$.

I claim that the measures $\mu_{N,t}$ converge weakly to a
measure $\mu_{0,t}$ on ${\mathbb R}^{k\nu}$ as $N\to\infty$. This
convergence implies the above inequality. This convergence can
be proved similarly to Proposition~5.1. Namely, we can write
\begin{eqnarray*}
&&\prod_{l=1}^\nu 
\frac{\exp\left\{i\frac{]t^{(l)}N[}N(x_1^{(l)}+\cdots+x_k^{(l)})\right\}-1} 
{N\left(\exp\left\{i\frac1N(x_1^{(l)}+\cdots+x_k^{(l)})\right\}-1
\right)} \\
&&\qquad =\frac1{N^\nu}\sum_{u\in B_N(t)}
\exp\left\{i\frac1N(u,x_1+\cdots+x_k)\right\},
\end{eqnarray*}
where the set $B_N(t)$ was defined in~(\ref{1.11}). With the help
of this formula the Fourier transform of the measure $\mu_{N,t}$
can be calculated in all points of the form $u=(u_1,\dots,u_k)$,
$u_s=\frac {p_s}N$, $p_s\in{\mathbb Z}^\nu$, $1\le s\le k$, 
This can be done similarly to the corresponding calculation in
Proposition~5.1. Then a good asymptotic formula can be proved
for this Fourier transform with the help of relation~(\ref{1.3}),
and this implies the above mentioned convergence. Here again the
method of proof in Proposition~5.1 is applied. I omit the details.

This implies that condition~(b) of Proposition~4A holds in our
model. The proof of the remaining conditions is much simpler.
Similarly to the proof of Theorem~3.2 it can be shown
with the help of Proposition~3.1 that the spectral measures
$G^{(N)}_{j,j'}$ satisfy the required convergence property.
Finally, it is not difficult to check that the functions
$h^N_{j_1,\dots,j_k}$ defined in~(\ref{2.25}) satisfy
condition~(a) of Proposition~4A. Theorem~3.4 is proved.

\medskip
Let me finally remark that a simple and natural
modification in the proof of  Theorem~3.4 shows that this
result also holds if the random variables $S_0(t)$
in it are defined for all $t\in[0,\infty)^\nu$, (in the
way as it is explained at the end of Section~3) and not
only for $t\in[0,1]^\nu$. 

\appendix
\section{On the background of the limit theorems of this paper.}
\label{A.}

In the example after formula~(\ref{1.3}) I constructed a vector
valued stationary random field with a spectral density
function in such a way that its covariance function satisfies
relation~(\ref{1.3}) with some appropriately defined functions
$a_{j,j'}\left(\frac p{|p|}\right)$ and $L(p)$. The spectral
density of this random field is close in some sense to the
spectral density of a vector valued generalized stationary
random field. Moreover, the spectral density of this
generalized random field has some homogeneity property.

I would like to point out that the spectral measures of all
stationary Gaussian random fields whose covariance function
satisfy~(\ref{1.3}) show a similar behavior. Indeed, take
the spectral measure $(G_{j,j'}(\cdot))$, $1\le j,j'\le d$,
of such a random field whose covariance function $r_{j,j'}(p)$,
$1\le j,j'\le d$, $p\in{\mathbb Z}^\nu$  satisfies
relation~(\ref{1.3}). Let us recall the results of
Proposition~3.1 about the properties of this spectral measure.

The elements, $G_{j,j'}^{(N)}(\cdot)$, defined as
$G^{(N)}_{j,j'}(A)=\frac{N^\alpha}{L(N)}G_{j,j'}\left(\frac AN\right)$,
$1\le j,j'\le d$, $N=1,2,\dots$, $A\subset {\mathbb R}^\nu$, of
the rescaled versions of the spectral measure $G=(G_{j,j'}(\cdot))$,
$1\le j,j'\le d$, of a random field, whose covariance matrix
satisfies~(\ref{1.3}) have a vague limit $G^{(0)}_{j,j'}(\cdot)$, when
$N\to\infty$. These vague limits have the homogeneity property
$G^{(0)}_{j,j'}(tA)=t^\alpha G^{(0)}_{j,j'}(A)$, $1\le j,j'\le d$,
for all $t>0$ and measurable, bounded sets $A\subset{\mathbb R}^\nu$.
Moreover, $(G^{(0)}_{j,j'}(\cdot))$, $1\le j,j'\le d$, is the spectral
measure of a generalized, stationary random field.

The above mentioned homogeneity property of the measures
$G^{(0)}_{j,j'}$ is important for us, because it enables us to
construct self-similar random fields, and in our limit theorems
self-similar random fields appear as the limit. Here I recall
the definition of self-similarity in a slightly more general
situation than in the main text. In this definition vector
valued random fields are considered. A vector valued random
field $S(t)=(S_1(t),\dots,S_m(t))$, $t\in[0,\infty)^\nu$, of
dimension~$m$ is called self-similar with parameter~$\beta$,
$\beta>0$, if $S(ut)\stackrel{\Delta}{=}u^\beta S(t)$ for all
$u>0$, where $\stackrel{\Delta}{=}$ means that the finite
dimensional distributions of the two random fields agree.

To understand how a vector valued self-similar random field
can be constructed with the help of the spectral measure
$(G_{j,j'}(\cdot))$ of a vector valued stationary generalized
random field whose elements have the homogeneity property
$G_{j,j'}(tA)=t^{\alpha}G_{j,j'}(A)$ with some $\alpha>0$ let us
first recall that the set of functions $\varphi$ for which
the random variable $Z_G(\varphi)$ of a generalized, vector
valued Gaussian random field with spectral measure $G$ is
defined can be enlarged. Indeed, let
$Z_G=(Z_{G,1},\dots,Z_{G,d})$ be a random spectral measure
corresponding to the spectral measure $(G_{j,j'}(\cdot))$.
In the original definition the elements of the vector
valued, generalized Gaussian stationary random field
corresponding to this random spectral measure are the
random vectors $(Z_{G,1}(\varphi),\cdots,Z_{G,d}(\varphi))$,
with $Z_{G,j}(\varphi)=\int \tilde\varphi(x)Z_{G,j}(\,dx)$,
$1\le j\le d$, where the function $\varphi(\cdot)$ is an
element of the Schwartz space ${\cal S}$, and
$\tilde\varphi$ denotes its Fourier transform. These
random integrals $Z_{G,j}(\varphi)$ can be defined for
a larger class of functions. They can be defined for
those real valued functions $\varphi(x)$, for which
$\int |\tilde\varphi(x)|^2\,G_{j,j}(\,dx)<\infty$
for all $1\le j\le d$. 

The spectral measure $(G_{j,j'}(\cdot))$, $1\le j,j'\le d$,
of those generalized random fields are important for us
for which the domain of arguments of the random variables
$Z_{G,j}(\varphi)$, $1\le j\le d$, can be extended with the
indicator functions of the rectangles
$[0,t]=\prod\limits_{s=1}^\nu[0,t_s]$ for all
$t=(t_1,\dots,t_\nu]\in{\mathbb R}^\nu$, i.e.
$$
\int |\widetilde{I_{[0,t]}}(x)|^2 G_{j,j}(\,dx)=
\int \left(\prod_{s=1}^\nu\frac{2(1-\cos (t_sx_s))}{x_s^2}\right)
G_{j,j}(\,dx)<\infty
$$
for all $1\le j\le d$. This inequality holds if
$G_{j,j}(tA)=t^{\alpha}G_{j,j}(A)$ with some $0<\alpha<2\nu$
for all $t>0$ and $1\le j\le d$.

Let us consider a spectral measure $(G_{j,j'}(\cdot))$,
$1\le j,j'\le d$, such that
$$
G_{j,j'}(tA)=t^{\alpha}G_{j,j'}(A) \textrm{  with some } 0<\alpha<2\nu
$$
for all $t>0$ and $1\le j,j'\le d$, and let
$(Z_{G,1}^\alpha,\dots,Z_{G.d}^\alpha)$
be a random spectral corresponding to this spectral measure. (Here I
put the homogeneity parameter $\alpha$ of the spectral measure
in the upper index of the elements of the random spectral measure.)
Consider for all pairs of vectors
$t^{(1)}=(t_1^{(1)},\dots,t_\nu^{(1)})\in{\mathbb R}^\nu$ and
$t^{(2)}=(t_1^{(2)},\dots,t_\nu^{(2)})\in{\mathbb R}^\nu$ such that
$t^{(1)}_s<t^{(2)}_s$ for all $1\le s\le\nu$ the rectangle
$[t^{(1)},t^{(2)}]=\prod\limits_{s=1}^\nu [t^{(1)}_s,t^{(2)}_s]$, and
define the random vectors
$$
Z^{\alpha}_G([t^{(1)},t^{(2)}])=
(Z_{G,1}^\alpha([t^{(1)},t^{(2)}]),\dots,Z_{G.d}^\alpha([t^{(1)},t^{(2)}]))
$$
with coordinates
\begin{eqnarray*}
Z_{G,j}^\alpha([t^{(1)},t^{(2)}])
&=& \int\widetilde{I_{[t^{(1)},t^{(2)}]}}(x)Z_{G,j}^\alpha(\,dx) \\
&=& \int\left(\prod_{s=1}^\nu
\frac{e^{i(t^{(2)}_sx_s}-e^{it^{(1)}_s)x_s}}{ix_s}\right)
Z_{G,j}^\alpha(\,dx),\quad 1\le j\le d,
\end{eqnarray*}
for all these rectangles.

Introduce the vectors $S_0(t)=Z^{\alpha}_G([0,t])$ for all
$t\in{\mathbb R}^\nu$ with positive coordinates, where $0$ denotes
the origin in ${\mathbb R}^\nu$, and $X_0(p)=Z^{\alpha}([p-1,p])$
for all $p\in{\mathbb Z}^\nu$, where $p-1=(p_1-1,\dots,p_\nu-1)$
for $p=(p_1,\dots,p_\nu)$. Then $S_0(\cdot)$ is a vector valued
self-similar random field with self-similarity parameter
$\nu-\frac\alpha2$, $X_0(p)$, $p\in{\mathbb Z}^\nu$, is a vector
valued stationary, Gaussian random field, and for a fixed vector
$p=(p_1,\dots p_\nu)\in{\mathbb Z}^\nu$ and all $N=1,2,\dots$
$$
\frac1{N^{\nu-\alpha/2}} \sum_{\substack{j=(j_1,\dots,j_\nu)\\
1\le j_s\le Np_s\textrm{ for all }1\le s\le\nu}} X_0(j)=  
\frac1{N^{\nu-\alpha/2}}S_0(Np)\stackrel{\Delta}{=}S_0(p).
$$

A similar relation holds also for the linear combinations of the
coordinates of the vector valued random field $X_0(p)$,
$p\in{\mathbb Z}^\nu$. This means that these random fields satisfy
the limit theorems of Theorems~3.2---3.4 for $k=1$. Let me remark
that the covariance function $r_{j,j'}(p)=EX_0(0)X_0(p))$,
$p\in{\mathbb Z}^\nu$, satisfies relation~(\ref{1.3}). Indeed, it
can be proved that
\begin{eqnarray*}
r_{j,j'}(p)&=&\int \widetilde{I_{[0,1]}}(x)\widetilde{I_{[p,p+1]}}(x)
G_{j,j'}(\,dx) \\
&=&\int e^{i(p,x)}\left(\prod_{s=1}^\nu\frac{2(1-\cos x_s)}{x_s^2}\right)
G_{j,j'}(\,dx) \\
&=&C_{j,j'}\left(\frac p{|p|}\right)|p|^{-\alpha}(1+o(1))
\end{eqnarray*}
with some function $C_{j,j'}\left(\frac p{|p|}\right)$
because of the homogeneity property of~$G_{j,j'}(\cdot)$.

Theorems 1.2---1.4 in the case $k=1$ state that the
corresponding limit theorems also hold for models which satisfy
Condition~(\ref{1.3}) with $0<\alpha<\nu$. The restriction of
the value of $\alpha$ to $0<\alpha<\nu$ instead of
$0<\alpha<2\nu$ in these results has a good reason.
The partial sums which were normalized in these theorems have
variances of order $N^{2\nu-\alpha}L(N)$. In the case $\alpha>\nu$
this means an exponent smaller than~$\nu$. So in the case
$\alpha>\nu$ we can get a limit theorem (for $k=1$) only in such
models where both positive and negative covariances appear, and
their effects compensate each other in a very special way.

\medskip
In the case $k>1$ a similar picture arises. Here again, we take
random spectral measures corresponding to such spectral measures
which have homogeneity property. We define the self-similar
random fields we are working with in this case by means of
$k$-fold Wiener--It\^o integrals with respect to random
spectral measures corresponding to them.

Given an integer $k\ge2$ let us take the spectral measure
$(G_{j,j'}(\cdot))$ of a generalized random field which has the
homogeneity property $G_{j,j'}(tA)=t^{\alpha}G_{j,j'}(A)$ for all
measurable sets with finite diameter, $t>0$, $1\le j,j'\le d$,
and a number $\alpha>0$ whose possible value will be given later.
Let us consider a random spectral measure $Z_G=(Z_{G,1},\dots,Z_{G,d})$
corresponding to this spectral measure, and define with its help
a random vector $Z_G([t^{(1)},t^{(2)}])$ for all rectangles
$[t^{(1)},t^{(2)}]$ introduced in the previous construction
with coordinates $Z_{G,j_1,\dots,j_k}([t^{(1)},t^{(2)}])$, 
where $(j_1,\dots,j_k)$  is a sequence with the property
$1\le j_s\le d$ for all $1\le s\le k$. These random variables are
defined as $k$-fold Wiener--It\^o integrals by the following
formula.
\begin{eqnarray*}
&&Z_{G,j_1,\dots,j_k}([t^{(1)},t^{(2)}])
=\int\widetilde{I_{[t^{(1)},t^{(2)}]}} (x^{(1)}+\dots+x^{(k)}) \\
&&\qquad\qquad\qquad\qquad\qquad\qquad\qquad\qquad
Z_{G,j_1}(\,dx^{(1)})\dots  Z_{G,j_k}(\,dx^{(k)})\\
&&\qquad=\int\left(\prod_{s=1}^\nu
\frac{e^{ it^{(2)}_s(x_s^{(1)}+\dots+x_s^{(k)})}
-e^{ it^{(1)}_s(x_s^{(1)}+\dots+x_s^{(k)})}}
{i\left(x_s^{(1)}+\dots+x_s^{(k)}\right)}\right) \\
&&\qquad\qquad\qquad\qquad\qquad\qquad\qquad\qquad
Z_{G,j_1}(\,dx^{(1)})\dots  Z_{G,j_k}(\,dx^{(k)}).
\end{eqnarray*}
(Actually, the value of this random integral depends only on the
multiplicity of the numbers $1,2,\dots,d$ in the sequence
$(j_1,\dots,j_k)$. The order of these numbers in the sequence
$j_1,\dots,j_k$ does not count.)

With the help of the above defined Wiener--It\^o integrals 
let us define, similarly to the case $k=1$, the random vector
$$
Z_G([t_1,t_2])=\{Z_{G,j_1,\dots,j_k}([t^{(1)},t^{(2)}])\colon\;
1\le j_s\le d \textrm{ for all }  1\le s\le l\}.
$$

Naturally, we have to choose the spectral measure $(G_{j,j'}(\cdot)$
in such a way that the above $k$-fold Wiener--it\^o integral be
meaningful. There is an integral which must be finite for the
existence of these random integrals. In the case of $k$-fold
Wiener--It\^o this relation holds if the homogeneity parameter
$\alpha$ of the underlying spectral measure satisfies the
inequality $0<\alpha<\frac\nu k$. The next calculation is an
estimation which implies the existence of the above
Wiener--it\^o integrals  with such a choice of~$\alpha$. I omit
the explanation why this calculation is correct, although
this is not self-evident. Actually the existence of the random
integrals I am considering here also follows from the results
of the main text.

In the case $\alpha<\frac\nu k$
\begin{eqnarray*}
&&\int\left(\prod_{s=1}^\nu
\frac{1-\cos(t^{(2)}_s-t^{(1)}_s)(x_s^{(1)}+\dots+x_s^{(k)})}
{\left(x_s^{(1)}+\dots+x_s^{(k)}\right)^2}\right) \\
&& \qquad\qquad\qquad         
G_{j_1, j_1}(\,dx^{(1)})\dots  G_{j_k,j_k}(\,dx^{(k)}) \\
&&\qquad \le C\int\left(\prod_{s=1}^\nu
\frac 1{\left(1+\left(x_s^{(1)}+\dots+x_s^{(k)}\right)^2\right)}\right) \\ 
&& \qquad\qquad\qquad         
|x^{(1)}|^{\alpha-\nu}\cdots |x^{(k)}|^{\alpha-\nu} \,dx^{(1)}\dots\,dx^{(k)
} \\ 
&&\qquad \le C'\prod_{s=1}^\nu\left(\int
\frac {|x_s^{(1)|}|^{-1+\alpha/\nu}\cdots|x_s^{(k)}|^{-1+\alpha/\nu}}
{\left(1+\left(x_s^{(1)}+\dots+x_s^{(k)}\right)^2\right)}
\,dx_s^{(1)}\dots\,dx_s^{(k)}\right) <\infty.
\end{eqnarray*}

Similarly to the case $k=1$ such a special limit theorem
can be presented which shows some similarity to the results of
Theorems~3.2---3.4. Define, with similar notation as in the case
$k=1$, the random vectors $S_0(t)=Z_G([0,t])$ for all
$t\in{\mathbb R}^\nu$ with positive coordinates, and
$X_0(p)=Z^{\alpha}([p-1,p])$ for all $p\in{\mathbb Z}^\nu$.
Then $S_0(\cdot)$ is a vector valued self-similar random field
with self-similarity parameter $\nu-\frac{k\alpha}2$,  $X_0(p)$,
$p\in{\mathbb Z}^\nu$, is a vector valued stationary, Gaussian
random field, and for a fixed vector
$p=(p_1,\dots p_\nu)\in{\mathbb Z}^\nu$, $p_s>0$,$\le \nu$,
and all $N=1,2,\dots$
$$
\frac1{N^{\nu-k\alpha/2}} \sum_{\substack{j=(j_1,\dots,j_\nu)\\
1\le j_s\le Np_s\textrm{ for all }1\le s\le\nu}} X_0(j)=  
\frac1{N^{\nu-k\alpha/2}}S_0(Np)\stackrel{\Delta}{=}S_0(p).
$$

This is a limit theorem, where the limit is the above constructed
self-similar random field $S_0(\cdot)$. Theorems~3.2--3.4 are limit
theorems with the same limit. They hold for such vector valued
Gausian stationary  random fields whose covariance matrices are
similar to the covariance matrix of a random field with that
spectral measure $(G_{j,j'}(\cdot)$, $1\le j,j'\le d$, which was
applied in the definition of the self-similar random field
$S_0(\cdot)$. 

\medskip
The goal of this Appendix was to explain the background of the
results in this paper. Here I concentrated on the explanation of
the definition of the self-similar random fields which appear as the
limit in our limit theorems. Their construction was based on the
theory of multiple Wiener--It\^o integrals for several dimensional
stationary Gaussian random fields, in particular on the properties
of the random spectral measures of generalized Gaussian random fields.

The proof of the results consisted of two steps. In the first
step the random sums we wanted to study were rewritten as a
sum of Wiener--It\^o integrals. This could be done with the
help of the multivariate version of It\^o's formula formulated
in Theorem~2.2 of~\cite{11}. Then limit theorems for sequences
of sums of Wiener--It\^o integrals had to be proved. This
could be done with the help of Proposition~4A. Naturally, the
Wiener--It\^o integrals which define the limit random variables
in Theorems~3.2--3.4 must exist. The condition
$0<\alpha<\frac{\nu}k$ in these theorems appear because of this
condition. They are to guarantee the existence of the $k$-fold
Wiener--It\^o integrals which define the limit random variables
in these results.

If the conditions of Theorems~3.2--3.4 hold, but with a
parameter $\alpha\ge\frac{\nu}k$, then the random integrals
defining the limit in these theorems do not exist. In such
cases the central limit theorem holds with the classical
normalization. This follows from the result of~\cite{3} in
the scalar valued and from its multivariate generalization
in Theorem~4 of~\cite{1} in the vector valued case. This
problem is discussed in Appendix~B of~\cite{11}.

\medskip
In this paper limit theorems were proved for non-linear
functionals of stationary Gaussian random fields. I try to give
a short overview about papers which deal with similar problems.
The scalar valued version of the results in this paper was
proved in~\cite{6}. M.~S.~Taqqu proved in paper~\cite{16}
similar results. Both papers contain non-central limit theorems
for sequences of non-linear functionals of (scalar valued)
stationary Gaussian random fields. Taqqu's result has no
multivariate version, and I do not know how such a result
can be proved.

Paper~\cite{2} contains a result which tells when the classical
central limit theorem holds for sequences of non-linear
functionals of stationary Gaussian random fields similar to
those considered in~\cite{6}. The book~\cite{12} generalizes
the result of this paper. A.~M.~Arcones wanted to generalize
both the central limit theorem of~\cite{2} and the non-central
limit theorem of~\cite{6} for non-linear functionals of vector
valued Gaussian random fields. He proved the central limit
theorem part of these results in Theorem~4 of his paper~\cite{1}.
He claimed to have also proved the multivariate version of the
non-central limit theorem in~\cite{6} in Theorem~6 of his paper.
But I found the proof (and also the formulation) of this result
problematic. The goal of the present paper was to formulate and
prove the right multivariate version of the result in~\cite{6}.

Let me also remark that H.~C.~Ho and T.~C.~Sun proved an
interesting result in~\cite{8}. They proved a result  which
can be considered as an interesting mixture of the central and
non-central limit theorem for non-linear functionals of stationary
Gaussian random fields. They considered a two-dimensional
vector valued stationary Gaussian random process $(X_m,Y_m)$,
$m\in{\mathbb Z}$, together with two non-linear functionals
which are of such type as the non-linear functionals applied
in~\cite{6}. They applied the first functionals for the
process~$X_m$, and the second functional for the process~$Y_m$.
In Theorem~1 of their paper they investigated the case when
the sequence of the linear functionals defined with the help
of the elements in the first coordinate satisfy a non-central,
and the corresponding sequence defined with the help of the
elements in the second coordinate satisfy the central limit
theorem. They proved under some additional conditions that the 
joint distributions of these sequences also have a limit, and
the two coordinates of this limit are independent. (The processes
$X_m$ and $Y_m$ are not independent.) Let me remark that the
usual proofs of the central and non-central limit theorem apply
different methods. This indicates that the proof of this result
in~\cite{8} demanded new ideas.

Donatas Surgailis proved results about similar problems, and they
are also worth mentioning. He proved, together with some
coauthors such results which can be considered as a generalization
of the central and non-central limit theorems proved for
non-linear functionals of stationary Gaussian random fields.
He proved limit theorems for non-linear functionals of a
new class of stationary stochastic processes which contains
non-Gaussian processes, too. He worked with scalar valued random
processes, but probably these results can be generalized also to
vector valued random processes.

Surgailis has several articles about this subject. 
I would mention his paper~\cite{16}, where he investigated
non-linear functionals of moving average processes which may
be non-Gaussian. He proved limit theorems for non-linear
functionals of such processes. These results are very similar
to those proved for Gaussian processes. The main point in the
proofs of this paper is that Surgailis considered the Appel
polynomials related to the moving average process he was
working with, and showed that the arguments applied in the
Gaussian case can be adapted to the problems investigated by
him if the Hermite polynomials are replaced by these Appel
polynomials. (In the case of Gaussian moving averages the
Appel polynomials are Hermite polynomials.)

\medskip
Finally I briefly mention a field of research where similar limit
theorems are proved with the help of essentially different arguments.
This research deserves special attention because of its importance.
This is the theory about the KPZ (Kardar--Parisi--Zhang) universality
classes and the problems related to them. Many important problems
can be studied with their help. On the other hand, the application
of this theory demands hard analysis. An overview about it can be
found in paper~\cite{4} together with a long list of literature.

\section{Proof of the corollary of Theorem~3.4.}
\label{B.}

{\it Proof of the corollary  of Theorem~3.4.}\/
We want to show that for all pairs $\varepsilon>0$ and $\eta>0$
there exists some $\delta=\delta(\varepsilon,\eta)>0$ and threshold
index $N_0=N_0(\varepsilon,\eta)$ such that for all $N\ge N_0$ the
inequality
\begin{equation}
P\left(\sup_{\substack{(s,t)\colon\;s,t\in[0,1]^\nu\\
\;\;\; |t-s|<\delta}}|S_N(t)-S_N(s)|>\varepsilon\right)<\eta \label{B1}
\end{equation}
holds for the random field $S_N(t)$, $t\in[0,1]^\nu$, defined
in~(\ref{1.11}) and~(\ref{1.12}).

Inequality~(\ref{B1}) means that the random fields $S_N(t)$,
$t\in[0,1]^\nu$, introduced in Theorem~3.4 satisfy besides the limit
theorem formulated in Theorem~3.4 also the tightness condition for
probability measures in the space of continuous functions
$C([0,1]^\nu,{\cal C})$. It can be seen that these two properties
together imply the desired functional limit theorem.

\medskip
First I show that relation~(\ref{B1}) can be replaced
with the following set of simpler inequalities.

Define for all $1\le j\le\nu$ and $\delta>0$ the  following
set $V(j,\delta)$ consisting of pairs
of vectors $(s,t)$, $s\in[0,1]^\nu$ and $t\in[0,1]^\nu$.
\begin{eqnarray*}
V (j,\delta)
\!\!&=&\!\! \{(s,t)\colon\;
s=(s_1,\dots,s_\nu),\,t=(s_1,\dots,s_{j-1},
t_j,s_{j+1},\dots,s_\nu) \colon\; \\
&& \quad 0\le s_l\le 1 \textrm{ for all }1\le l\le \nu,
\textrm{ and } s_j\le t_j\le\min(1,s_j+\delta)\}.
\end{eqnarray*}
Then
\begin{equation}
P\left(\sup_{(s,t)\in V(j,\delta)}
|S_N(t)-S_N(s)|>\varepsilon\right)<\eta \quad \textrm{for all }
1\le j\le\nu  \label{B2}
\end{equation}
if $\delta\le\delta(\varepsilon,\eta)$ and
$N\ge N_0(\varepsilon,\eta)$
with some $\delta(\varepsilon,\eta)>0$ and
$N_0(\varepsilon,\eta)$.

\medskip
To see the possibility of such a reduction let us first observe
that inequality~(\ref{B1}) follows from its following formally
weaker version.

For all $\varepsilon>0$ and $\eta>0$ there exists
some $\delta=\delta(\varepsilon,\eta)>0$ and
$N_0=N_0(\varepsilon,\eta)$ such that
\begin{equation}
P\left(\sup_
{\substack{(s,t)\colon\;0\le s_j<t_j\le1, \\
\qquad \textrm{ for all }1\le j\le\nu, \;\; |t-s|<\delta}}
|S_N(t)-S_N(s)|>\varepsilon\right)<\eta \label{B3}
\end{equation}
if $N\ge N_0$.

Indeed, for a pair of vectors $(s,t)$, $s\in[0,1]^\nu$,
$t\in[0,1]^\nu$, define the vector
$$
s^*=s^*(s,t)=(\min(s_1,t_1),\dots,\min(s_\nu,t_\nu)),
$$
and consider the pairs $(s^*,s)$ and $(s^*,t)$. Let us apply
relation~(\ref{B3}) with parameters $\delta$ and $N_0$ corresponding
to the parameters $\frac\varepsilon2$ and $\frac\eta2$. If the pair $(s,t)$
satisfies the conditions appearing in the supremum of~(\ref{B1}) with
these parameters, then the pairs $(s^*,s)$ and $(s^*,t)$ satisfy
the conditions in the supremum of~(\ref{B3}) with the same
parameter~$\delta$. Also the relation
$|S_N(t)-S_N(s)|\le |S_N(s)-S_N(s^*)|+|S_N(t)-S_N(s^*)|$ holds. Hence
inequality~(\ref{B3}) with the above chosen $\delta$ and $N_0$
implies that
\begin{eqnarray*}
 && \!\!\!\!\!\!\!\!\!\!\!\!
 \sup_{\substack{(s,t)\colon\; s,t\in[0,1]^\nu \\ |t-s|<\delta }} |S_N(t)-S_N(s)|
\le\sup_ {\substack{(s^*,s)\colon\;0\le s^*_j<s_j\le1, \\
 \quad \textrm{ for all }1\le j\le\nu, \;\; |s-s^*|<\delta}}
|S_N(s)-S_N(s^*)|  \\
&&\qquad  +\sup_{\substack{(s^*,t)\colon\;0\le s^*_j<t_j\le1, \\
 \quad \textrm{ for all }1\le j\le\nu,\;\; |t-s^*|<\delta}}
|S_N(t)-S_N(s^*)|\le\varepsilon
\end{eqnarray*}
with probability more than $1-\eta$ if $N>N_0$. This means that
relation~(\ref{B3}) implies relation~(\ref{B1}).

Relation~(\ref{B2}) can be reduced to reation~(\ref{B3}) in a
similar way. To do this let us first define for a pair of vectors
$s=(s_1,\dots,s_\nu)\in[0,1]^\nu$, $t=(t_1,\dots,t_\nu)\in[0,1]^\nu$ and
number $1\le j\le\nu$ the vector
$s(j)=s(j,s,t)=(t_1,\dots,t_{j-1},s_j,\dots,s_\nu)$
(for $j=1$ $s(1)=(s_1,\dots,s_\nu)=s$), and consider the pairs of vectors
$(s(j),s(j+1))$, $1\le j\le \nu-1$. Observe that
$(s(j),s(j+1))=(s(j,s,t),s(j+1,s,t))\in V(j,\delta)$ if the pair
$(s,t)$ satisfies the relations
$|t-s|<\delta$ and $0\le s_j<t_j\le1$ for all $1\le j\le\nu$. Let us
choose a $\delta>0$ and $N_0$ in such a way that inequality~(\ref{B2})
holds with parameters $\frac\varepsilon\nu$ and $\frac\eta\nu$ with this
number~$\delta$ and $N\ge N_0$. Let us take those pairs of vectors
$(s,t)$ which satisfy the conditions imposed in the supremum of
formula~(\ref{B3}) with this number $\delta$. We have seen that
\begin{eqnarray*}
&&\{(s(j,s,t),s(j+1,s,t))\colon\; 0\le s_l<t_l\le1, \\
&& \qquad\qquad  \textrm{ for all } 0\le l\le\nu,
|t-s|\le\delta\}\subset V(j,\delta).
\end{eqnarray*}
The identity
$S_N(t)-S_N(s)=\sum\limits_{j=1}^{\nu-1}[S_N(s(j+1,s,t))-S_N(s(j,s,t))]$
also holds. These relations imply that with our choice of $\delta$
$$
\sup_{\substack{(s,t)\colon\;0\le s_j<t_j\le1, \\
\qquad \textrm{ for all }1\le j\le\nu, \;\; |t-s|<\delta}}
|S_N(t)-S_N(s)|\le\sum_{j=1}^{\nu-1} \sup_{(s,t)\in V(j,\delta)}|S_N(t)-S_n(s)|
\le\varepsilon
$$
with probability more than $1-\eta$ if $N>N_0$. This means that
(\ref{B2}) implies (\ref{B3}).

Next I present an inequality with the help of the random variables
\hfill\break
$H(X_1(p),\dots,X_d(p))$ instead of $S_N(t)$ which implies 
inequality~(\ref{B2}). For this goal I introduce the following
notations.

Let us define the rectangle $D_N(r,s)$ 
for all pairs of vectors $r=(r_1,\dots,r_\nu)$ and
$s=(s_1,\dots,s_\nu)$ with integer coordinates such that
$0\le r_j< s_j\le N$ for all $1\le j\le \nu$ by the formula
$$
D_N(r,s)=\{p=(p_1,\dots,p_\nu)\colon\; p\in{\mathbb Z}^{\nu},
\quad r_j< p_j\le s_j \textrm{ for all } 1\le j\le \nu\},
$$
and introduce for all $\delta>0$ and $1\le j\le \nu$ a set
${\cal D}_N(\delta,j)$ consisting of the above defined rectangles
$D_N(r,s)$ with some additional properties. We define
\begin{eqnarray*}
  {\cal D}_N(\delta,j)&=&\{D_N(r,s)\colon\; r_l=0
  \textrm { for } l\neq j, \;\; 0\le r_j\le N, \\
&&\quad 0<s_l\le N \textrm { for all } 1\le l\le \nu, \textrm{ and }
0<s_j-r_j\le \delta N\}.
\end{eqnarray*}
Inequality~(\ref{B2}) follows from the relation
\begin{eqnarray}
&&P\left(\sup_{D_N(r,s)\in {\cal D}_N(\delta,j)}\left|
\frac{\sum\limits_{p\in D_N(r,s)} \
H(X_1(p),\dots,X_d(p))}{N^{\nu-k\alpha/2}L(N)^{k/2}} \right|
>\varepsilon\right) \le  \eta \nonumber \\
&&\qquad\qquad\qquad\qquad\qquad \textrm{ for all } 1\le j\le\nu \label{B4}
\end{eqnarray}
if $\delta\le\delta(\varepsilon,\eta)$ with some
$\delta(\varepsilon,\eta)>0$. Here, and also in the remaining part
of the proof $H(x_1,\dots,x_d)$ is a sum of the form
$$
H(x_1,\dots,x_d)=H^{(0)}(x_1,\dots,x_d)+H^{(1)}(x_1,\dots,x_d),
$$
with functions $H^{(0)}(\cdot)$ and $H^{(1)}(\cdot)$ defined
in formulas~(\ref{1.5}) and~(\ref{1.10}).

In formulas~(\ref{B4}) and~(\ref{B2}) very similar expressions
are estimated. The main difference between them is that
in~(\ref{B2}) random variables of the form $|S_N(t)-S_N(s)|$
are considered with arguments $s,t\in[0,1]^\nu$, while
in~(\ref{B4}) random variables of the form
$\left|S_N\left(\frac sN\right)-S_N\left(\frac rN\right)\right|$
with arguments $\frac sN$ and $\frac rN$, where $s$ and $r$ are
vectors with integer coordinates with values between~0 and~$N$.
This is a sort of discretization, and in the reduction
of~(\ref{B4}) to~(\ref{B2}) it has to be shown that this
discretization has a negligible effect in the case of a large
sample size~$N$.

This can be seen with the help of the following observation. If $N$
is large, then because of the definition of the random field
$S_N(\cdot)$ for all $t\in[0,1]^\nu$ there exists a vector
$r=(r_1,\dots,r_\nu)$ with integer coordinates $r_j$, $0\le r_j \le N$,
such that $S_N(t)=S_N\left(\frac rN\right)$, and $t$ and $\frac rN$
are very close to each other.

\medskip
Inequality~(\ref{B4}) will be proved by means of a good estimate
on the tail distribution of the random variables
$\sum\limits_{p\in D_N(r,s)} H(X_1(p),\dots,X_d(p))$ for the rectangles
$D_N(r,s)$. These expressions will be estimated by means of an
argument similar to the proof of Theorem~3.3. To do this let us
first remark that Lemma~1 of~\cite{1} implies the following result,
too. The inequality in Lemma~5A holds also in the case when the
function~$H^{(1)}(x_1,\dots,x_d)$ is replaced in it by
$H(x_1,\dots,x_d)$, and the coefficient $\psi^{k+1}$ in the
upper bound is replaced by $\psi^k$.

This modified version of Lemma~5A yields that there exists a threshold
index $n_0$ and some constant $C>0$ such that if the parameters $p$ and
$q$ of two elements $X(p)=(X_1(p),\dots,X_d(p))$
and $X(q)=(X_1(q),\dots,X_d(q))$ of our random field satisfy the
inequality $|p-q|\ge n_0$, then
$$
|EH(X_1(p),\dots X_d(p))H(X_1(q),\dots,X_d(q))|
\le C|p-q|^{-k\alpha}L(|p-q|)^{k}.
$$ 

Let us observe that for such pairs $p$ and $q$ also the inequality
$$
|p-q|^{-k\alpha}L(|p-q|)^k\le C\prod_{\substack{j\colon\;1\le j\le \nu\\
 \;\;\;\;\; p_j-q_j\neq0}} |p_j-q_j|^{-k\alpha/\nu}
(L(|p_j-q_j|)^{k/\nu}+I_{\{|p_j-q_j|<D\}})
$$
holds with some $C>0$ and $D>0$, where $I_{\{x<D\}}$ denotes the
indicator function of the set $\{x\colon\; x<D\}$. Hence the
previous estimate has the consequence
\begin{eqnarray*}
&&|EH(X_1(p),\dots X_d(p))H(X_1(q),\dots,X_d(q))| \\
&&\qquad \le C\prod_{\substack{j\colon\; 1\le j\le \nu\\
 \;\;\;\;\;\; p_j-q_j\neq0}}
|p_j-q_j|^{-k\alpha/\nu} (L(|p_j-q_j|)^{k/\nu}+I_{\{|p_j-q_j|\le D\}})
\end{eqnarray*}
if $|p-q|\ge n_0$. This inequality is more appropriate for us than
the previous one.

On the other hand, the inequality
\begin{eqnarray*}
&&|EH(X_1(p),\dots, X_d(p))H(X_1(q),\dots,X_d(q))|\\
&&\qquad \qquad \le E H^2(X_1(0),\dots X_d(0))\le C_1
\end{eqnarray*}
also holds for all $p,q\in{\mathbb Z}^\nu$ with some $C_1<\infty$
because of the Schwarz inequality and relation~(\ref{1.10}).

The last two inequalities imply that for any rectangular
$D_N(r,s)\subset B_N$ and $p\in D_N(r,s)$
\begin{eqnarray}
&&\left|EH(X_1(p),\dots,X_d(p))
\left(\sum_{q\in D_N(r,s)} H(X_1(q),\dots,X_d(q))\right)\right|
\label{B5} \\
&&\qquad\qquad\qquad \le C_2\prod_{j=1}^\nu
\left(1+(s_j-r_j)^{1-k\alpha/\nu}L(s_j-r_j)^{k/\nu}\right) \notag
\end{eqnarray}
with an appropriate constant~$C_2$. Indeed, these inequalities
imply that 
\begin{eqnarray*}
&&\left|EH(X_1(p),\dots,X_d(p))
\left(\sum_{q\in D_N(r,s)} H(X_1(q),\dots,X_d(q))\right)\right| \\
&&\qquad\qquad \le C\prod_{j=1}^\nu\left(1+2\sum_{q_j=1}^{s_j-r_j}
q_j^{-k\alpha/\nu}(L(q_j)^{k/\nu}+I_{\{q_j\le D\}})\right) 
\end{eqnarray*}
and $\sum\limits_{q_j=1}^{s_j-r_j} q_j^{-k\alpha/\nu}
(L(q_j)^{k/\nu}+I_{\{q_j<\le D\}})
\le C'(1+ (s_j-r_j)^{1-k\alpha/\nu} L(s_j-r_j)^{k/\nu})$ with some
$C'>0$, since $k\alpha/\nu<1$ by the conditions of Theorem~3.4.
These relations imply~(\ref{B5}).

By summing up inequality~(\ref{B5}) for all $p\in D_N(r,s)$, and
applying an appropriate normalization we get that
\begin{eqnarray}
&&\frac1{N^{2\nu-k\alpha} L(N)^k}
E\left[\sum_{p\in D_N(r,s)} H(X_1(p),\dots,X_d(p)\right]^2 \label{B6} \\
&&\qquad\qquad \le C_2\prod_{j=1}^\nu
\left(\frac{(s_j-r_j)+(s_j-r_j)^{2-k\alpha/\nu} L(s_j-r_j)^{k/\nu}}
{N^{2-k\alpha/\nu}L(N)^{k/\nu}}\right). \notag
\end{eqnarray}

I claim that if $1\le s_j-r_j\le N$ for some $1\le j\le\nu$ and
$\eta>0$ chosen so small that $\beta=\frac{k\alpha}\nu+\eta<1$, then
\begin{equation}
\frac{(s_j-r_j)+(s_j-r_j)^{2-k\alpha/\nu} L(s_j-r_j)^{k/\nu}}
 {N^{2-k\alpha/\nu}L(N)^{k/\nu}}\le
C\left(\frac{s_j-r_j}N\right)^{2-k\alpha/\nu-\eta} \label{B7}
\end{equation}
with some $C=C(\eta)>0$.

Indeed,
\begin{eqnarray*}
\frac{(s_j-r_j)}{N^{2-k\alpha/\nu}L(N)^{k/\nu}}&=&
\left(\frac{s_j-r_j}N\right)^{2-k\alpha/\nu-\eta}
 (s_j-r_j)^{k\alpha/\nu+\eta-1}\frac{N^{-\eta}}{L(N)^{k/\nu}} \\
 &\le& C\left(\frac{s_j-r_j}N\right)^{2-k\alpha/\nu-\eta}
\end{eqnarray*}
if $\eta>0$ is chosen so small that $\frac{k\alpha}\nu+\eta<1$,
and
\begin{eqnarray*}
&&\frac{(s_j-r_j)^{2-k\alpha/\nu} L(s_j-r_j)^{k/\nu}}
 {N^{2-k\alpha/\nu}L(N)^{k/\nu}}=
\left(\frac{s_j-r_j}N\right)^{2-k\alpha/\nu-\eta} \\
&&\qquad \left(\frac{s_j-r_j}N\right)^\eta
\left(\frac{L(s_j-r_j)}{L(N)}\right)^{k/\nu}
\le C\left(\frac{s_j-r_j}N\right)^{2-k\alpha/\nu-\eta}.
\end{eqnarray*}
These two inequalities imply~(\ref{B7}).

Let us choose a sufficiently large number $D>0$ (whose
value does not depend on $N$ and $\delta$), and introduce the
quantity $d_N(p)=\frac DN$ for all $1\le p\le N$. With such
a notation we can get the following inequality with the help
of relations~(\ref{B6}) and~(\ref{B7}).

Take some rectangle $D_N(r,s)\subset B_N$. Then we have
for any $\lambda>0$
\begin{eqnarray}
&& P\left(\frac
{\left|\sum\limits_{p\in D_N(r,s)} H(X_1(p),\dots,X_d(p))\right|}
{N^{\nu-k\alpha/2} L(N)^{k/2}}>\lambda \right) \label{B8} \\ 
&&\qquad \le \frac1{\lambda^2}
\frac{ E\left[\sum\limits_{p\in D_N(r,s)} H(X_1(p),\dots,X_d(p))\right]^2}
{N^{2\nu-k\alpha} L(N)^k}   \notag \\
&&\qquad \le\frac{C_3}{\lambda^2}\prod_{j=1}^\nu
\left(\frac{s_j-r_j}N\right)^{2-k\alpha/\nu-\eta}
\le\frac1{\lambda^2}\prod_{j=1}^\nu
\left(\sum_{p_j=r_j+1}^{s_j} d_N(p_j)\right)^{2-\beta} \notag
\end{eqnarray}
with $\beta=\frac{k\alpha}\nu+\eta<1$.

With the help of formula~(\ref{B8}) one can get such a maximum-type
inequality which implies formula~(\ref{B4}).

In the case $\nu=1$ Theorem~10.2 of Billingsley's book~\cite{2}
can be applied. In this case this result together with
formula~(\ref{B8}) imply that
\begin{eqnarray*}
&& P\left(\sup_{(u,v)\colon\;r<u<v\le s}
\left|\frac {\sum\limits_{u<p\le v} H(X_1(p),\dots,X_d(p))}
{N^{\nu-k\alpha/2} L(N)^{k/2}}\right|>\lambda \right) \\
 &&\qquad \le \frac K{\lambda^2}   
\left(\sum_{r<p\le s} d_N(p)\right)^{2-\beta}
=\frac {D^{2-\beta}Kn}{\lambda^2}\left(\frac{s-r}N\right)^{2-\beta} 
\end{eqnarray*}
for any pairs $0\le r<s\le N$ with some $K>0$. In particular,
$$
P\left(\sup_{(u,v)\colon\;r<u<v\le r+\delta}
\left|\frac {\sum\limits_{u<\frac pN\le v} H(X_1(p),\dots,X_d(p))}
{N^{\nu-k\alpha/2} L(N)^{k/2}}\right|>\lambda \right) 
\le \frac {K'}{\lambda^2}\delta^{2-\beta}   
$$
for any interval $[r,r+\delta]\subset[0,1]$ with $\delta>0$,

Since the exponent of $\delta$ in the last inequality equals
$2-\beta>1$ it is not difficult to see that this relation implies
inequality~(\ref{B4}) in the case $\nu=1$.

Indeed, we get~(\ref{B4}) by applying this inequality with
the choice $\lambda=\varepsilon$ for the intervals
$[k\delta,(k+2)\delta]$ for all $0\le k<\frac1\delta$ with a
sufficiently small $\delta>0$. Then inequality~(\ref{B4})
implies inequality~(\ref{B1}), too.

\medskip
There is a multivariate version of the inequality cited from
Billingsley's book~\cite{2} also in the case $\nu>1$ which,
together with formula~(\ref{B8}) imply inequality~(\ref{B4})
in the general case. This inequality implies for any $\nu\ge1$
that if inequality~(\ref{B8}) holds for all
rectangles $D_N(r,s)$, then
\begin{eqnarray}
 && P\left(\sup_{\substack{D_N(u,v)\colon\; \\
 D_N(u,v)\subset D_N(r,s)}}
\left|\frac{\sum\limits_{p\in D_N(u,v)} H(X_1(p),\dots,X_d(p))}
{N^{\nu-k\alpha/2} L(N)^{k/2}}\right|>\lambda \right) \label{B9} \\
&&\quad \le\frac K{\lambda^2}
\prod_{j=1}^\nu\left(\sum_{r_j <p_j\le s_j} d_N(p_j)\right)^{2-\beta} 
=\frac {D^{(2-\beta)\nu} K}{\lambda^2}
\prod_{j=1}^\nu\left(\frac{s_j-r_j}N\right)^{2-\beta} \notag
\end{eqnarray}
with some $K>0$, $D>0$ and $\beta=\frac{k\alpha}\nu+\eta<1$.
(Here we are working with the previously defined $d_N(p)=\frac DN$.)

Indeed, although I did not find this result in the literature
there is such a generalized version of the inequality quoted from
Billingsley's book which states that if inequality~(\ref{B8})
holds, then it implies inequality~(\ref{B9}), too. This can be
proved for instance by means of induction with respect to the
dimension~$\nu$ by exploiting that this result holds for $\nu=1$.
In the proof we have to exploit that the upper bound
in~(\ref{B8}) has a special product form.

\medskip
Let us fix some parameter $1\le j\le\nu$, a number $0<\delta\le 1$
an integer $0\le r\le N$, and define with their help, similarly to
the definition of ${\cal D}_N(\delta,j)$ the set of
of rectangles
\begin{eqnarray*}
 {\cal D}_N(\delta,j,r)&=&\{D_N(u,v)\colon\; u_l=0,\textrm{ and }
0<v_l\le N \textrm { for } l\neq j, \\
&&\qquad\qquad \textrm{ and } \;  r\le u_j<v_j\le r+N\delta\}.
\end{eqnarray*}

The following inequality is a special case of ~(\ref{B9}).
\begin{eqnarray*}
&&P\left(\sup_{D(u,v)\in{\cal D}_N(\delta,j,r)}
\left|\frac{\sum\limits_{p\in D_N(u,v)} H(X_1(p),\dots,X_d(p))}
{N^{\nu-k\alpha/2} L(N)^{k/2}}\right|>\lambda \right) \\
&& \qquad\qquad\qquad\qquad
\le\frac{K'}{\lambda^2}\left(\frac{s_j-r_j}N\right)^{2-\beta}. 
\end{eqnarray*}

Inequality~(\ref{B4}), hence inequality~(\ref{B1}) can be
proved with the help of the last inequality in the same way as
it was done for $\nu=1$. Actually, in that proof a special case
of this inequality was applied. Since, as it was mentioned at
the start of the proof relation~(\ref{B1}) together with
Theorem~3.4 imply the desired weak convergence the corollary
is proved. 

\medskip\noindent
Acknowledgements. The author got support from the Hungarian 
Foundation NKFI--EPR No. K-125569.

\end{document}